\documentclass{article}

\PassOptionsToPackage{square, numbers}{natbib}

\usepackage[final]{neurips_2020}

\usepackage[utf8]{inputenc} 
\usepackage[T1]{fontenc}    
\usepackage{hyperref}       
\usepackage{url}            
\usepackage{booktabs}       
\usepackage{amsfonts}       
\usepackage{nicefrac}       
\usepackage{microtype}      
\usepackage{amsmath}
\usepackage{macros}
\usepackage{mathtools}
\usepackage{algorithm}
\usepackage{algorithmic}

\usepackage{graphicx}
\usepackage{subfigure}
\newenvironment{proof}{\paragraph{Proof:}}{\hfill$\square$}
\usepackage{xcolor}

\usepackage{enumitem}
\usepackage{changes}
\definechangesauthor[name={Per cusse}, color=red]{imp}
\usepackage{bbm}
\usepackage{caption}
\usepackage{multicol}
\usepackage{multirow}

\allowdisplaybreaks[4]

\title{Biased Stochastic First-Order Methods for Conditional Stochastic Optimization and Applications in Meta Learning
}

\author{
  	Yifan Hu\thanks{The first two authors have equal contribution.}\\
	UIUC\\
	\texttt{yifanhu3@illinois.edu}
	\And
	Siqi Zhang\footnotemark[1]\\
	UIUC\\
	\texttt{siqiz4@illinois.edu}
	\AND
	Xin Chen\\
	UIUC\\
	\texttt{xinchen@illinois.edu}
	\And
	Niao He\\
    UIUC \& ETH Zurich\\
    \texttt{niao.he@inf.ethz.ch}
}

\begin{document}

\maketitle

\begin{abstract}
Conditional stochastic optimization covers a variety of applications ranging from invariant learning and causal inference to meta-learning. However, constructing unbiased gradient estimators for such problems is challenging due to the composition structure. As an alternative, we propose a biased stochastic gradient descent (BSGD) algorithm and study the bias-variance tradeoff under different structural assumptions. We establish the sample complexities of BSGD for strongly convex, convex, and weakly convex objectives under smooth and non-smooth conditions. Our lower bound analysis shows that the sample complexities of BSGD cannot be improved for general convex objectives and nonconvex objectives except for smooth nonconvex objectives with Lipschitz continuous gradient estimator.  For this special setting, we propose an accelerated algorithm called biased SpiderBoost (BSpiderBoost) that matches the lower bound complexity.  We further conduct numerical experiments on invariant logistic regression and model-agnostic meta-learning to illustrate the performance of BSGD and BSpiderBoost.
\end{abstract}

\section{Introduction}
We study a class of optimization problems, called \emph{conditional stochastic optimization (CSO)}:
\begin{equation}
\label{pro:ori}
    \min_{x\in\mathcal{X}} F(x) := \EE_\xi f_\xi(\EE_{\eta|\xi} g_\eta(x,\xi)),
\end{equation} 
where $\mathcal{X}\subseteq \RR^d$, $g_\eta(\cdot,\xi):\RR^d\to\RR^k$  is a vector-valued function dependent on both random vectors $\xi$ and $\eta$,  $f_\xi(\cdot):\RR^k\to\RR$ depends  on the random vector $\xi$, and the inner expectation is taken with respect to the conditional distribution of $\eta|\xi$. Throughout, we assume access to samples from the distribution $P(\xi)$ and the conditional distribution $P(\eta|\xi)$.  

CSO includes the classical stochastic optimization as a special case when $g_\eta(x,\xi)=x$  but is much more general. It has been recently utilized to solve a variety of applications in machine learning, ranging from the policy evaluation and control in reinforcement learning~\citep{pmlr-v54-dai17a,pmlr-v80-dai18c,nachum2020}, the optimal control in linearly-solvable Markov decision process~\citep{pmlr-v54-dai17a}, to instrumental variable regression in causal inference~\citep{singh2019kernel,muandet2019dual}. 

One common challenge with these applications is that in the extreme case, a few or only one sample is available from the conditional distribution of $\eta|\xi$ for each given $\xi$. To deal with this limitation, a primal-dual stochastic approximation algorithm was proposed to solve a min-max reformulation of CSO using the kernel embedding techniques~\citep{pmlr-v54-dai17a}. However, this approach requires convexity of $f$ and linearity of $g$, which are not satisfied by general applications when neural networks are involved. 

On the other hand, for many other applications, e.g., those arising from invariant learning and meta-learning, we do have access to multiple samples from the conditional distribution. Take the model-agnostic meta-learning (MAML)~\citep{finn2017model} as an example. 
MAML learns a meta-initialization parameter using metadata from similar learning tasks such that taking one or multiple gradient steps on a small training data would generalize well on a new task. It can be framed into the following CSO problem:
\begin{equation}
\label{pro:maml}
    \min_{w}\; \EE_{i\sim p, a\sim D_{query}^i} l_i\Big(\EE_{b\sim D_{support}^i}\big(w-\alpha\nabla l_i(w,b)\big),a\Big),
\end{equation}
where $p$ represents the distribution of different tasks,  $D_{support}^i$ and $D_{query}^i$ correspond to support (training) data and query (testing) data of the task $i$, $l_i(\cdot,D^i)$ is the loss function on data $D^i$ from task $i$, and $\alpha$ is a fixed meta step size. Setting $\xi = (i,a)$ and $\eta = b$, \eqref{pro:maml} is clearly a special case of CSO for which multiple samples can be drawn from the conditional distribution of $P(\eta|\xi)$. Since the loss function generally involves neural networks, the resulting CSO is often  nonconvex. Thus, the previous primal-dual algorithm and kernel embedding techniques developed in~\citep{pmlr-v54-dai17a} no longer apply.

In this paper, we focus on the general CSO problem where multiple samples from the conditional distribution $P(\eta|\xi)$ are available, and the objective is not necessarily in the compositional form of a convex loss $f_\xi(\cdot)$ and a linear mapping $g_\eta(\cdot,\xi)$. Recently, \citet{hu2019sample} studied the generalization error bound and sample complexity of empirical risk minimization (ERM), a.k.a., sample average approximation (SAA) for general CSO:
$$
\min_{x\in\mathcal{X}} \frac{1}{n}\sum_{i=1}^n f_{\xi_i}\Big(\frac{1}{m}\sum_{j=1}^m g_{\eta_{ij}}(x,\xi_i)\Big),
$$
where  $\{\xi_i\}_{i=1}^n$ are i.i.d. samples from $\PP(\xi)$ and $\{\eta_{ij}\}_{j=1}^m$ are i.i.d. samples from $\PP(\eta|\xi_i)$. They assumed that the global optimal solution to ERM can be computed without specifying how.   Differently, here we aim at developing efficient stochastic gradient-based methods that directly solve the CSO problem~\eqref{pro:ori} and find either a global optimal solution in the convex setting or a stationary point in the nonconvex setting, respectively.

Due to the composition structure of the CSO objective in~\eqref{pro:ori}, constructing unbiased gradient estimators is not possible in general. Instead, we  leverage a mini-batch of conditional samples to construct the gradient estimator with controllable bias and propose a family of biased first-order methods, including (1) the biased stochastic gradient descent (BSGD) algorithm for general convex and nonconvex CSO objectives and (2) the biased SpiderBoost (BSpiderBoost) algorithm, designed for nonconvex smooth CSO objectives. Note that BSpiderBoost is inspired by the variance reduced method for nonconvex smooth stochastic optimization in \citep{fang2018spider, wang2018spiderboost}.

\subsection{Our contributions}
Our main results are summarized in Table \ref{tab:result}. Our contributions are three-fold:\\[-0.5cm]

\begin{itemize}[leftmargin=2em]
    \item We establish the first sample complexity results of BSGD and BSpiderBoost in the context of CSO. Since the bias of BSGD comes from estimating the conditional expectation rather than from a given stochastic oracle, the sample complexity closely depends on the smoothness conditions of the outer function, which is distinct from traditional SGD results. For convex problem, to achieve an $\eps$-optimal solution, 
    the sample complexity of BSGD improves from $\cO(\eps^{-4})$ to  $\tilde \cO(\eps^{-3})$ when either $f_\xi$ is smooth or $F$ is strongly convex and further improves to  $\tilde \cO(\eps^{-2})$ when both conditions hold, where $\tilde \cO(\cdot)$ represents the bound with hidden logarithmic factors.  For weakly convex CSO problems, BSGD requires a total sample complexity of $ \cO(\eps^{-8})$ to achieve an $\eps$-stationary point, and of $\cO(\eps^{-6})$ when $f_\xi$ is smooth. If we further assume that both $f_\xi$ and $g_\eta$ are Lipschitz continuous and Lipschitz smooth, the biased gradient estimator is Lipschitz continuous, and then the sample complexity can be improved to $\cO(\eps^{-5})$ by BSpiderBoost.

\item {\color{black} We analyze the lower bounds on the minimax error of first-order algorithms using specific biased oracles for CSO objectives. With the upper bounds results, BSGD is optimal for strongly convex, convex, and weakly convex CSO objectives, and BSpiderBoost is optimal for the nonconvex smooth CSO problems under the additional oracle assumption that the gradient estimator returned by the oracle is Lipschitz continuous.
 }

\item When applied to MAML,  BSGD converges to a stationary point under simple deterministic stepsize rules and appropriate inner mini-batch sizes. In contrast, the commonly used first-order MAML algorithm~\citep{finn2017model} ignores the Hessian information and is not guaranteed to converge even when a large inner mini-batch size is used. For smooth MAML, compared with the algorithm recently introduced in \citep{fallah2019convergence}, BSGD  without requiring stochastic stepsizes and mini-batches of outer samples at each iteration, thus is more practical. Leveraging the variance reduction technique, BSpiderBoost attains the best-known sample complexity for MAML, to our best knowledge.  Numerically, we further demonstrate that BSGD and BSpiderBoost achieve superior performance for MAML.  
\end{itemize}

\begin{table*}[t]
    {
    \caption{{\color{black}Sample Complexity for CSO }}
    \renewcommand\arraystretch{1.2}
    \begin{center}
    \begin{tabular}{ccccccccc}
        \toprule[1.5pt]
        
        \multirow{3}{*}{\textbf{Algorithm}}
        & \multicolumn{8}{c}{\textbf{Assumptions}}
        \\
        \cline{2-9}
        
        {}
        & \emph{$\hat F$}
        & \emph{SC}
        & \emph{SC}
        & \emph{Convex}
        & \emph{Convex}
        & \emph{WC}
        & \emph{WC}
        & \emph{Smooth}
        \\
        \cline{2-9}
        
        {}
        & \emph{$f_\xi$}
        & \emph{Smooth}
        & \emph{Lipschitz}
        & \emph{Smooth}
        & \emph{Lipschitz}
        & \emph{Smooth}
        & \emph{Lipschitz}
        & \emph{Smooth}
        \\
        \hline
        
        \multicolumn{2}{c}{ERM~\citep{hu2019sample}*}
        & \emph{$\cO(\eps^{-2})$}
        & \emph{$\cO(\eps^{-3})$}
        & \emph{$\tilde \cO(d\eps^{-3})$}
        & \emph{$\tilde \cO(d\eps^{-4})$}
        & -
        & -
        & -
        \\
        \hline
        
        \multicolumn{2}{c}{BSGD}
        & \emph{$\tilde \cO(\eps^{-2})$}
        & \emph{$\tilde \cO(\eps^{-3})$}
        & \emph{$ \cO(\eps^{-3})$}
        & \emph{$ \cO(\eps^{-4})$}
        & \emph{$ \cO(\eps^{-6})$}
        & \emph{$ \cO(\eps^{-8})$}
        & \emph{$ \cO(\eps^{-6})$}
        \\
        \hline
        
        \multicolumn{2}{c}{BSpiderBoost}
        & \emph{-}
        & \emph{-}
        & \emph{-}
        & \emph{-}
        & \emph{-}
        & \emph{-}
        & \emph{$\cO(\eps^{-5})$}
        \\
        \hline
        
        \multicolumn{2}{c}{Lower Bound }
        & \emph{$ \cO(\eps^{-2})$}
        & \emph{$ \cO(\eps^{-3})$}
        & \emph{$ \cO(\eps^{-3})$}
        & \emph{$ \cO(\eps^{-4})$}
        & \emph{$ \cO(\eps^{-6})$}
        & \emph{$ \cO(\eps^{-8})$}
        & \emph{$ \cO(\eps^{-5})$}
        \\
        \hline
        
        \multicolumn{9}{c}{\textbf{Goal}: find $\eps$-optimal solution for convex $F$ and $\eps$-stationary point for weakly convex $F$.} \\ 
        \multicolumn{9}{c}{$\hat F$ is defined in \eqref{eq:empirical}.
        SC: strongly convex; WC: weakly convex; Lipschitz = Lipschitz continuous.} \\
        \multicolumn{9}{c}{* SAA  requires further solving the empirical risk minimization.} \\
        \bottomrule[1.5pt]
    \end{tabular}  
    \vskip -0.2in
    \label{tab:result}
    \end{center}
}
\end{table*}

\subsection{Related work} 
\paragraph{Nested expectation optimization (NEO)} \citep{wang2016accelerating,wang2017stochastic,ghadimi2018single,yang2019multilevel,chen2020solving,zhang2019multi} deals with problems in the form of:
$\min_{x \in \mathcal{X}}\;  \mathbf{f}\circ \mathbf{g}(x):=\EE_{\xi}\big[f_\xi\big({\EE_{\eta}[g_\eta(x)]}\big)\big]$,  where  $\mathbf{f}(u):=\EE_{\xi}[f_\xi(u)]$, $\mathbf{g}(x):=\EE_{\eta}[g_\eta(x)]$. A key assumption is that even when $\xi$ and $\eta$ are dependent, there exists a \emph{deterministic function} $\mathbf{g}(x)$, independent of $\xi$, which does not hold for general CSO. Hence their algorithms  and analysis cannot extend to CSO.  The best known sample complexity for smooth strongly convex NEO is $\cO(\eps^{-1.25})$~\citep{wang2016accelerating}; and for nonconvex smooth NEO objective is $\cO(\eps^{-4})$~\citep{ghadimi2018single} and $\cO(\eps^{-3})$~\citep{zhang2019multi} if $\xi$ and $\eta$ are independent using variance reduction technique. 

\paragraph{Nested expectation estimation} Estimating nested expectations in the form of $
\EE[ H(\EE(\eta|\xi))]$ has been extensively studied in the statistics and simulation communities. \citep{hong2009estimating, gordy2010nested, hong2017kernel} considered nested Monte Carlo estimator, when $H$ is a general non-linear function. \citep{bujok2015multilevel,giles2015multilevel,giles2019multilevel} considered Multilevel Monte Carlo (MLMC) method~\citep{giles2008multilevel} when $H$ has special structure. Note that this line of work purely focuses on estimation, whereas we deal with optimization, which is more challenging.

\paragraph{Biased gradient methods} \citep{hu2016bandit, ajalloeian2020analysis, karimi2019non, hu2020analysis, chen2018stochastic} analyzed the non-asymptotic convergence of general biased gradient methods.  These papers assume that the bias in gradient estimator comes from certain black-box oracles or an additive non-zero mean noise. Differently in our problem, the bias directly comes from estimating the nested expectation and can be controlled by the sampling strategy.

\paragraph{Notations}
$\Pi_\mathcal{X}$ denotes the projection operator, i.e., $\Pi_\mathcal{X}(x) := \argmin_{z\in\mathcal{X}}\|z-x\|_2^2$. $\tilde \cO(\cdot)$ represents the order hiding logarithmic factors. A function $f(\cdot):\mathbb{R}^k\rightarrow\mathbb{R}$  is $L$-Lipschitz continuous on $\mathcal{X}$ if $|f(x)-f(y)|\leq L\|x-y\|_2$ holds for any $x,y\in\mathcal{X}$. A function $f(\cdot)$ is $S$-Lipschitz smooth on $\mathcal{X}$ if  $f(x)-f(y)-\nabla f(y)^\top (x-y)\leq \frac{S}{2}\|x-y\|_2^2$ holds for any $x,y\in\mathcal{X}$. 
A function $f(\cdot)$ is $\mu$-convex on $\mathcal{X}$ if for any $x,y\in\mathcal{X}$, 
$
    f(x)-f(y)-\nabla f(y)^\top(x-y)\geq \frac{\mu}{2}\|x-y\|_2^2.
    $
Note that $\mu >0$, $\mu=0$, and $\mu<0$ correspond to $f$ being strongly convex, convex, and weakly convex, respectively. Lastly, we denote  $x^*\in\argmin_{x\in\mathcal{X}}F(x)$ as an optimal solution to the problem of interest.  For an abuse of notation, we use $\nabla$ to denote the Jacobian matrix, (sub)gradient vector, and derivative for simplicity.

\section{Biased Stochastic First-Order Methods}
\label{section:BSGD}
For simplicity, throughout, we assume that $\cX\subseteq\RR^d$ is closed  and convex, and the random functions $f_\xi(\cdot)$ and $g_\eta(\cdot,\xi)$ are continuously differentiable. Based on the special composition structure of CSO and the chain rule, under mild conditions, the gradient of $F(x)$ in \eqref{pro:ori} is given by
\begin{equation*}
\nabla F(x) = \EE_\xi\Big[(\EE_{\eta|\xi}\nabla g_{\eta}(x,\xi)])^\top   \nabla f_\xi(\EE_{\eta|\xi}g_{\eta}(x,\xi))\Big].  
\end{equation*}
Constructing an unbiased stochastic estimator of the gradient can be costly and even impossible.  Instead, we consider a biased estimator of $\nabla F(x)$ using one sample $\xi$ and $m$ i.i.d. samples $\{\eta_{j}\}_{j=1}^{m}$ from the conditional distribution of $P(\eta|\xi)$ in the following form:
\begin{equation}
\label{eq:nabla_empirical}
\nabla \hat{F}(x; \xi,\{\eta_{j}\}_{j=1}^{m}
)
:=\Big(\frac{1}{m}{\sum}_{j=1}^m\nabla g_{\eta_j}(x,\xi)\Big)^\top \nabla f_\xi\Big(\frac{1}{m}{\sum}_{j=1}^m g_{\eta_j}(x,\xi)\Big).    
\end{equation}
Note that $\nabla \hat{F}(x; \xi,\{\eta_{j}\}_{j=1}^{m})$ is the gradient of an empirical objective such that
\begin{equation}
\label{eq:empirical}
\hat F(x;\xi,\{\eta_{j}\}_{j=1}^{m}) :=f_{\xi}\Big(\frac{1}{m}{\sum}_{j=1}^m g_{\eta_{j}}(x,{\xi})\Big).
\end{equation}
Based on this biased gradient estimator, we propose BSGD, which is formally described in Algorithm \ref{alg:A}. When using fixed inner mini-batch sizes $m_t=m$, BSGD can be viewed as performing SGD updates on the surrogate objective $\EE_{\{\xi,\{\eta_{j}\}_{j=1}^{m}\}}\hat F(x;\xi,\{\eta_{j}\}_{j=1}^{m})$. Inspired by the recent success of variance-reduced methods for nonconvex stochastic optimization~\citep{nguyen2017sarah,fang2018spider, wang2018spiderboost}, we further introduce an accelerated algorithm BSpiderBoost, which is formally described in Algorithm \ref{alg:BSpiderBoost}. BSpiderBoost divides updates into "epoch": at the beginning of the epoch, it will initialize the gradient estimator with $ N_1 $ outer samples of $ \xi $; then in later iterations in the epoch, the estimator will be updated with gradient information in current iteration generated with $ N_2 $ outer samples and the information from the last iteration. Compared to the classical SVRG method \citep{johnson2013accelerating}, this framework keeps utilizing the latest information for updates, which can generate more accurate gradient estimations.

\begin{algorithm}[t]
\caption{Biased Stochastic Gradient Descent (BSGD)}
\label{alg:A}
\begin{algorithmic}[1]
\REQUIRE Number of iterations $T$, inner mini-batch size $\{m_t\}_{t=1}^{T}$, initial point $x_1$, stepsize $\{\gamma_t\}_{t=1}^{T}$
\FOR{$t=1$ to $T-1$ \do} 
\STATE Sample $\xi_t$ from distribution $\PP(\xi)$, and $m_t$ i.i.d samples $\{\eta_{tj}\}_{j=1}^{m_t}$ from distribution $P(\eta|\xi_t)$.
\STATE Compute 
$\nabla \hat{F}(x_t; \xi_t,\{\eta_{tj}\}_{j=1}^{m_t})$ according to \eqref{eq:nabla_empirical}.
\STATE Update
        $$
        x_{t+1} = \Pi_\mathcal{X}\big(x_t-\gamma_{t} \nabla \hat{F}(x_t; \xi_t,\{\eta_{tj}\}_{j=1}^{m_t})\big).
        $$ 
\ENDFOR
\end{algorithmic}
\end{algorithm}

\begin{algorithm}[ht]
	\caption{Biased SpiderBoost (BSpiderBoost)}
	\label{alg:BSpiderBoost}
	\begin{algorithmic}[1]
		\REQUIRE  Number of iterations $ T $, inner batch size $m$, stepsize $\gamma$, epoch length $ q$, mini-batch sizes $ N_1, N_2,$
		\FOR{$t=0$ to $T$ \do} 
		\IF {$\mod(t,q)=0$}
		\STATE Generate $ N_1 $ samples of $\{\xi_1,\cdots,\xi_{N_1}\}$ 
		\STATE Generate $m$ i.i.d samples $\{\eta_{ij}\}_{i=1}^{m}$ from  $\mathbb{P}(\eta|\xi_i)$ for each $ \xi_i\in\{\xi_1,\cdots,\xi_{N_1}\} $.
		\STATE Compute $ v_t=\frac{1}{N_1}\sum_{i=1}^{N_1} \nabla \hat{F}(x_t; \xi_i,\{\eta_{ij}\}_{j=1}^{m})$ 
		\ELSE
		\STATE Generate $ N_2 $ samples of $\{\xi_1,\cdots,\xi_{N_2}\}$
		\STATE Generate $m$ i.i.d samples $\{\eta_{ij}\}_{i=1}^{m}$ from  $\mathbb{P}(\eta|\xi_i)$ for each $ \xi_i\in\{\xi_1,\cdots,\xi_{N_2}\} $.
		\STATE Compute 
		{\small
			\begin{equation}
			\label{eq:BSpB_Grad_Est}
    			v_t=\frac{1}{N_2}\sum_{i=1}^{N_2}\nabla \hat{F}(x_t; \xi_i,\{\eta_{ij}\}_{j=1}^{m})-\frac{1}{N_2}\sum_{i=1}^{N_2}\nabla \hat{F}(x_{t-1}; \xi_i,\{\eta_{ij}\}_{j=1}^{m})+v_{t-1}
			\end{equation}
		}
		\ENDIF
		\STATE Update $ x_{t+1}=x_t-\gamma v_t $
		\ENDFOR
		\ENSURE $ x_S $ which is uniformly randomly selected from $\{x_t\}_{t=1}^T$.
	\end{algorithmic}
\end{algorithm}

Before presenting the main results, we make one observation that the bias of the function value estimator $\hat F$,  
induced by the composition structure, depends on the smoothness condition of the outer function $f_\xi$:
{\color{black}for $L_f$-Lipschitz continuous $f_\xi$, 
$$
    \EE_{\xi, Y} [f_\xi(Y)-f_\xi(\EE_{Y|\xi} Y)] \leq L_f \EE_{\xi, Y|\xi} \|Y-\EE_{Y|\xi} Y\|_2;
$$ 
for $S_f$-Lipschitz smooth $f_\xi$,  
$$
    \EE_{\xi, Y} [f_\xi(Y)-f_\xi(\EE_{Y|\xi} Y)] \leq \frac{S_f}{2}\EE_{\xi, Y|\xi} \|Y-\EE_{Y|\xi} Y\|_2^2,  
$$
where $Y$ is a random variable. } To characterize the estimation error of $\hat F$, we make the following assumption.
\begin{assume} We assume that
\label{ass:general}
$\sigma_g^2 := {\sup}_{\xi,x\in\mathcal{X}} \EE_{\eta|\xi} \norm{ g_{\eta}(x,\xi)-\EE_{\eta|\xi} g_\eta(x,\xi)}_2^2<+\infty$;
\end{assume}
Assumption \ref{ass:general} indicates that the random vector $g_\eta$ has bounded variance.
Define 
\begin{equation}
\label{def:bias}
\Delta_f(m) = \left\{
\begin{array}{lr}
L_f\sigma_{g}/\sqrt{m},  \text{ if } f_\xi \text{ is } L_f\text{-Lipschitz continuous}, &\\
S_f\sigma_{g}^2/2m,  \text{ if } f_\xi \text{ is } S_f\text{-Lipschitz smooth}. \\
\end{array}
\right.
\end{equation}
The following lemmas characterize the estimation errors of $\hat F$ and $\nabla \hat F$. 

\begin{lm}[\citep{hu2019sample}]
\label{lm:bias} 
Under Assumption \ref{ass:general}, for a sample $\xi$ and $m$ i.i.d. samples $\{\eta_{j}\}_{j=1}^m$ from the conditional distribution $P(\eta|\xi)$, and any $x\in\cX$ that is independent of $\xi$ and $\{\eta_{j}\}_{j=1}^m$, we have
\beq{}
\big|\EE_{\{\xi,\{\eta_{j}\}_{j=1}^m \}}
\hat{F}(x; \xi,\{\eta_{j}\}_{j=1}^{m})   -   F(x)\big| \leq \Delta_f(m).
\eeq
\end{lm}
This implies that, to control the estimation bias up to $\eps$,   a  number of $m=\cO(\eps^{-2})$ samples is needed for Lipschitz continuous $f_\xi$ whereas  $m=\cO(\eps^{-1})$ is needed for Lipschitz smooth $f_\xi$.
{\color{black}\begin{lm}
\label{lm:bias_gradient}
Under Assumption \ref{ass:general}, if additionally assuming $f_\xi$ is $S_f$-Lipschitz smooth, $g_\eta$ is $L_g$-Lipschitz continuous, it holds that
\begin{equation}
\|\EE \nabla \hat{F}(x; \xi,\{\eta_{j}\}_{j=1}^{m}
    ) -\nabla F(x)\|_2^2\leq \frac{S_f^2L_g^2\sigma_g^2}{m}.
\end{equation}
\end{lm}}

\section{Convergence Analysis} 

In this section, we provide the non-asymptotic convergence analysis of BSGD and BSpiderBoost. 
Our result illustrates how the smoothness condition of the outer function $f_\xi$ influences the inner sample complexity and the total sample complexity.  This is quite different from the traditional SGD analysis for convex stochastic optimization, where the smoothness condition does not influence the complexity in terms of dependence on $\eps$~\citep{bubeck2015convex,nemirovski2009robust}. Before showing the convergence, we impose an assumption on the convexity.
  
\begin{assume}
\label{ass:convex}
$\hat{F}(x; \xi,\{\eta_{j}\}_{j=1}^{m})$ is $\mu$-convex for any $m,\xi,\{\eta_{j}\}_{j=1}^{m}$. 
\end{assume}

{Strong convexity, namely when $\mu>0$, can often be achieved by adding $\ell_2$-regularization to convex objectives. Convexity, namely when $\mu=0$, holds when  (i) $f_\xi$ is convex and $g_\eta$ is linear; (ii) $f_\xi$ and $g_\eta$ are convex and $f_\xi$ is non-decreasing. Weak convexity, namely when $\mu<0$, holds when  (i) $F$ is Lipschitz smooth, which holds if (i) both $f_\xi$ and $g_\eta$ are Lipschitz continuous and smooth  and (ii) $f_\xi$ is convex and $g_\eta$ is Lipschitz smooth~\citep{davis2019stochastic}.}
Note that weak convexity is commonly used in nonconvex optimization literature~\citep{paquette2017catalyst,chen2018universal,davis2019stochastic, zhang2018convergence}. Beyond weak convexity, little is known on the complexity of first-order algorithms except for some special functions, e.g., TAME functions \citep{davis2020stochastic} and difference of convex functions~\citep{pang2017computing}. Lastly, we point out that weak convexity is satisfied by various objectives used in machine learning, e.g., MAML discussed in this paper. Note that under mild conditions, Assumption \ref{ass:convex} implies that $F(x)$ is $\mu$-convex. 

\paragraph{Global convergence of BSGD for strongly convex objectives.} We have the following result:

\begin{thm}
\label{thm:stronglyconvex}
Under Assumption \ref{ass:general} and Assumption \ref{ass:convex} with $\mu>0$, if $\hat F$ is $S_F$-Lipschitz smooth, there exists a constant Set $\gamma_t  = \frac{1}{\mu(t+c)}$ with $c =\max\{ 4S_F^2/\mu^2-1,0\}$, the output $\hat x_T=\frac{1}{T}\sum_{t=1}^{T}x_t$ of BSGD
satisfies:
\begin{equation}
\EE[F(\hat{x}_T)-F(x^*)]\leq \frac{2\EE\|\nabla \hat F(x^*)\|_2^2(\log(T)+1)+ S_F^2/4\|x_1-x^*\|_2^2}{T\mu}+\frac{4}{T}\sum_{t=1}^T\Delta_f(m_t).
\end{equation}
\end{thm}
Hence, to achieve $\eps$-optimality, the number of iterations $T$ should be  at least $\tilde \cO(\eps^{-1})$, which aligns with the performance of SGD for strongly convex objectives~\citep{shamir2013stochastic}. 
For strongly-convex CSO with Lipschitz continuous $f_\xi$, recall the definition of $\Delta_f(m_t)$ in \eqref{def:bias}, using a fixed mini-batch size $m_t=\cO(\eps^{-2})$ or  time varying batch sizes $m_t=\cO(t^2)$ would be sufficient to obtain $\epsilon$-optimality. For strongly-convex CSO with Lipschitz smooth $f_\xi$,  it suffices to set $m_t=\cO(\eps^{-1})$ or  $m_t= \cO(t)$. Respectively, the total sample complexities are $\tilde \cO(\eps^{-3})$ and $\tilde \cO(\eps^{-2})$ under these two settings. 

\paragraph{Global convergence of BSGD for convex objectives}
We make an additional assumption about the second moment of the gradient estimator.
\begin{assume}
\label{ass:bounded_gradient}
There exists $M>0$ such that $
\EE \big[\|\nabla \hat{F}(x; \xi,\{\eta_{j}\}_{j=1}^{m}
)\|_2^2\mid x\big] \leq M^2
$
for any $x$.
\end{assume}
Note that Assumption \ref{ass:bounded_gradient} is common in the literature for analyzing SGD when the objective is non-strongly-convex or nonsmooth, e.g., when $f_\xi$ and $g_\eta$ are $L_f$- and $L_g$- Lipschitz continuous, $M = L_fL_g$. See e.g., \citep{nemirovski2009robust, bubeck2015convex,davis2019stochastic}.
\begin{thm}
\label{thm:convex}
Under Assumptions \ref{ass:general} Assumption \ref{ass:convex} with $\mu=0$, and Assumption \ref{ass:bounded_gradient}, with stepsizes $\gamma_t = c/\sqrt{T}$ for a positive constant $c$,   the output $\hat x_T=\frac{1}{T}\sum_{t=1}^{T}x_t$ of BSGD
satisfies 
\begin{equation*}
\EE [F(\hat x_T) -F(x^*)]\leq \frac{M^2c^2+\|x_1-x^*\|_2^2 }{2c\sqrt{T}}+\frac{2}{T}\sum_{t=1}^T\Delta_f(m_t).
\end{equation*}
\end{thm}
Comparing to Theorem \ref{thm:stronglyconvex}, without strong convexity condition, the iteration complexity increases  from $\cO(\eps^{-1})$ to $\cO(\eps^{-2})$. The total sample complexity for convex CSO is $\cO(\eps^{-4})$ for Lipschitz continuous $f_\xi$ with $m_t=\cO(\eps^{-2})$ or $m_t=\cO(t)$ and $\cO(\eps^{-3})$ for Lipschitz smooth $f_\xi$ with $m_t=\cO(\eps^{-1})$ or $m_t=\cO(\sqrt{t})$.

\paragraph{Stationary convergence of BSGD for general nonconvex objectives}
When $F(x)$ is nonconvex and possibly nonsmooth, we first introduce the notion of  convergence using Moreau envelope $F_{\lambda}(x)$ ($ \lambda>0 $) of  function $F(x)$ and its corresponding minimizer: 
\begin{equation*}
	\begin{split}
	F_{\lambda}(x):=\min_{z\in\mathcal{X}}\big\{F(z)+\frac{1}{2\lambda}||z-x||_2^2\big\}, \,\,
	\mathrm{prox}_{\lambda F}(x):=\underset{z\in\mathcal{X}}{\mathrm{argmin}}\big\{F(z)+\frac{1}{2\lambda}||z-x||_2^2\big\}.
	\end{split}
\end{equation*}
Based on Moreau envelope, we define the gradient mapping:
$
    \mathcal{G}_{\lambda F}(x):=\frac{1}{\lambda}||\mathrm{prox}_{\lambda F}(x)-x||_2.
$ 
We say $x$ is an $\eps$-stationary point of $F$ if  $\EE [\mathcal{G}_{\lambda F}(x)]\leq \eps$.
This convergence criterion is commonly used in nonconvex optimization literature \citep{beck17firstorder, drusvyatskiy2017proximal}. We have the following result.

\begin{thm}
\label{thm:wc_convergence}
	Under Assumption \ref{ass:general}, Assumption \ref{ass:convex} with $\mu <0$, and Assumption \ref{ass:bounded_gradient},  with  stepsizes $\gamma_t=c/\sqrt{T}$ for  a positive constant $c$, the output, $\hat x_R$, selected uniformly randomly from  $\{x_1,\cdots,x_{T}\}$, satisfies 
	\begin{equation*}\label{eq:nonconvex}
	\mathbb{E}\big[\mathcal{G}^2_{\frac{1}{2|\mu|}F}(\hat x_R)\big]
	\leq
	\frac{2F_{1/(2|\mu|)}(x_1)-2F(x^*)+2|\mu| M^2c^2}{{c\sqrt{T}}}
	+\frac{8|\mu|}{T}\sum_{t=1}^T\Delta_f(m_t).
	\end{equation*} 
\end{thm}
To the best of our knowledge, this is the first non-asymptotic convergence guarantee for CSO in the nonconvex setting. 
Specifically, for nonconvex CSO with Lipschitz continuous $f_\xi$,  setting batch sizes $m_t = \cO(\eps^{-4})$ or  $m_t = \cO(t)$ yields a total sample complexity of $\cO(\eps^{-8})$; for nonconvex CSO with Lipschitz smooth $f_\xi$,  using $m_t = \cO(\eps^{-2})$ or $m_t = \cO(\sqrt{t})$ achieves the total sample complexity of $\cO(\eps^{-6})$. Note that $\cO(\eps^{-6})$ sample complexity also holds when the CSO objective is additionally smooth. The analysis only requires little modification and is omitted.  

{\color{black}\paragraph{Stationary convergence of BSpiderBoost for nonconvex smooth objectives}
We now analyze the stationary convergence of BSpiderBoost for CSO problem with $\mathcal{X}=\RR^d$ and Lipschitz smooth $ F(x)$. We say $x$ is an $\eps$-stationary point of $F$ if $\EE \|\nabla F(x)\|_2\leq \eps$. Before proceeding,  we make the following assumption:
\begin{assume}
\label{ass:Lipschitz}
$f_\xi(\cdot)$ is $L_f$-Lipschitz continuous and $S_f$-Lipschitz smooth for any $\xi$. $g_\eta(\cdot,\xi)$ are  $L_g$-Lipschitz continuous and $S_g$-Lipschitz smooth for any $\xi$ and $\eta$.
\end{assume}
This assumption ensures that $F$ and $\hat F$ are $S_F$-Lipschitz smooth with $S_F = S_g L_f + S_f L_g^2$.

\begin{thm}[Convergence of BSpiderBoost]
	\label{thm:BSpiderBoost_CSO_revise}
	Under Assumptions \ref{ass:general} and \ref{ass:Lipschitz}, $\mathcal{X}=\RR^d$, and $ \Delta\coloneqq F(x_0)-F^* < \infty$, consider the following setup: $ T=\lceil 8\Delta\beta^{-1} \epsilon^{-2} \rceil $, $ q=\lfloor\sqrt{N_1}\rfloor $, $N_2=\big\lceil2\sqrt{N_1}\big\rceil$, $
	\gamma_t\equiv\gamma=\frac{1}{2S_F}$, and
	$$
	N_1=\Bigg\lceil\Big(3+\frac{3}{2\beta S_F}+\frac{3}{16\beta S_F}\Big)\frac{4L_f^2L_g^2}{\epsilon^2}\Bigg\rceil,\quad 
	m=\Bigg\lceil\Big(3+\frac{3}{2\beta S_F}+\frac{3}{16\beta S_F}\Big)\frac{4L_g^2S_f^2\sigma_g^2}{\epsilon^2}\Bigg\rceil,
	$$
	where
	$$
	\beta
	\coloneqq
	\frac{\gamma(1-S_F\gamma)}{2}-\frac{\gamma^3S_F^2q}{N_2}
	\geq
	\frac{1}{16S_F}>0.
	$$
	The output of BSpiderBoost, $ x_S $, which is randomly drawn from $\{x_1,...,x_T\}$, attains
$
	\mathbb{E}\|\nabla F(x_S)\|_2\leq \epsilon.
$
    Correspondingly, the sample complexity of BSpiderBoost is $ \mathcal{O}(\epsilon^{-5}) $.
\end{thm}

\begin{remark}
Recall that the objective of MAML \eqref{pro:maml} is a special case of CSO. If $\nabla l_i$ and $\nabla^2 l_i$ is Lipschitz continuous, then the objective is smooth and the outer function is smooth. Thus BSGD converges to an $\eps$-stationary point of \eqref{pro:maml} with sample complexity $ \cO(\eps^{-6})$ and BSpiderBoost converges with sample complexity $ \cO(\eps^{-5})$.
\end{remark}
}

\section{Lower Bounds for Conditional Stochastic Optimization}
\label{section:lower_bounds}
In this section, we show that the sample complexity of BSGD  for (strongly) convex and general nonconvex CSO objectives cannot be improved without further assumptions. The sample complexity achieved by BSpiderBoost also cannot be improved for nonconvex smooth CSO objectives.   The analysis uses the well-known oracle model, which consists of three components: a function class of interest $\mathcal{F}$, an algorithm class $\mathcal{A}$, and an oracle class $\Phi$. Specifically, we consider the biased stochastic first-order oracle class for CSO denoted as $\Phi_{m}$, where $m$ is the fixed number of conditional samples used. We also consider $\Phi_m^c$, a subset of $\Phi_m$ such that any oracle $\phi$

\begin{defn}[Biased first-order oracle for CSO]
For a query at point $x$ of CSO objective $F$ given by an algorithm, an oracle $\phi\in\Phi_m$ with a parameter $\sigma^2$ takes a sample $\zeta$ from its associated distribution $P(\zeta)$, and returns to the algorithm $\phi(x,F) = (h(x,\zeta),G(x,\zeta))$ 
such that
\begin{equation*}
\begin{split}
\EE h(x,\zeta) = \EE \hat F(x;\xi,\{\eta_{j}\}_{j=1}^m);  \quad  
 \EE G(x,\zeta)= \nabla \EE h(x,\zeta); \quad \EE\| G(x,\zeta)- \EE G(x,\zeta)\|_2^2 \leq \sigma^2.
\end{split}
\end{equation*}
In addition, we define the oracle class $\Phi_m^c$ such that $\Phi_m^c\subset \Phi_m$ and any oracle $\phi\in\Phi_m^c$ will return to the algorithm a Lipschitz continuous gradient estimator $G(\cdot,\zeta)$ for any $\zeta$.
\end{defn}

\textbf{Function class} we use $\mathcal{F}_\mathrm{CSO}$ to denote the CSO function class of interest. Specifically, we use $\cF_\mathrm{CSO}$ with superscript $+, 0, -$ to denote strongly convex, convex, nonconvex function class. 

\textbf{Randomized Algorithm}  A randomized algorithm class $\mathcal{A}$ contains algorithms $A$ such that $A$ maps the oracle output and a random seed $r$ to the next query point
$
x_{t+1}^A(\phi)=A(r,G(x_t^A(\phi), \zeta)).
$

\textbf{Updating Procedure} Suppose an algorithm $A\in\mathcal{A}$ is applied to minimize a function $F\in\mathcal{F}$ using oracle $\phi\in\Phi$. The updating procedure is such that at iteration $0$, the algorithm starts with some initialization point $x_0$. At iteration $t$, the algorithm $A$ queries the oracle $\phi$ about the information about $F$ on $x_t^A(\phi)$. The oracle $\phi$ will return some (noisy) information $\phi(x_t^A(\phi),F)$  back to the algorithm. Then the algorithm would base on all previous information returned by the oracle to generate the next query point $x_{t+1}^A(\phi)$. 

For a fixed number of iteration $T$, we define the minimax error as:
\begin{equation}
\label{eq:minimax}
\begin{split}
    & \Delta_T^*(\mathcal{A}, \mathcal{F},\Phi)
    := 
    \inf_{A \in \mathcal{A}} \sup_{\phi\in\Phi}\sup_{F \in \mathcal{F}}   \Delta_T(A, F,\phi):=\EE F(x_T^A(\phi))-\min_{x\in\mathcal{X}}F(x);\\
    & \Delta_T^{*g}(\mathcal{A}, \mathcal{F},\Phi)
    := 
    \inf_{A \in \mathcal{A}} \sup_{\phi\in\Phi}\sup_{F \in \mathcal{F}}   \Delta_T^{g}(A,F,\phi):=\EE\| \nabla F(x_T^A(\phi))\|_2^2,
\end{split}
\end{equation}
where the expectation is taken with respect to the randomness in algorithm $A$ and oracle $\phi$. $\Delta_T^*$ is used to capture the global optimality for convex function classes and  $\Delta_T^{*g}$ is used to capture the stationarity of the output for nonconvex function classes. If $\Delta_T^*\geq \eps$, it implies that for any algorithm $A$, there exists a `hard' function $F$ and an oracle $\phi$ such that the expected optimization error incurred by $A$   is at least $\eps$.

\begin{thm}
\label{thm:lower_bound_CSO}
For CSO problem, the minimax error satisfies that
\begin{itemize}[leftmargin=2em]
    \item[(i)] when $f_\xi$ is Lipschitz continuous,
    \begin{equation*}
    \begin{split}
    &\Delta_T^*(\mathcal{A}, \mathcal{F}_\mathrm{CSO}^+,\Phi_m)\geq \cO(m^{-1/2} + \sigma^2 T^{-1}); \   
    \Delta_T^*(\mathcal{A}, \mathcal{F}_\mathrm{CSO}^0,\Phi_m)\geq \cO( m^{-1/2} + \sigma T^{-1/2});\\
    &\Delta_T^{*g}(\mathcal{A}, \mathcal{F}_\mathrm{CSO}^-,\Phi_m)\geq \cO(m^{-1/2} + \sigma T^{-1/2}).
    \end{split}
    \end{equation*}
    \item[(ii)] when $f_\xi$ is Lipschitz smooth,
    \begin{equation*}
    \begin{split}
    & \Delta_T^*(\mathcal{A}, \mathcal{F}_\mathrm{CSO}^+,\Phi_m)\geq \cO(m^{-1} + \sigma^2 T^{-1}); \ 
    \Delta_T^*(\mathcal{A}, \mathcal{F}_\mathrm{CSO}^0,\Phi_m)\geq \cO(m^{-1} + \sigma T^{-1/2});\\
    & \Delta_T^{*g}(\mathcal{A}, \mathcal{F}_\mathrm{CSO}^-,\Phi_m)\geq \cO(m^{-1} +  \sigma T^{-1/2}).
    \end{split}
    \end{equation*}
    \item[(iii)] when the gradient estimator is Lipschitz continuous and $f_\xi$ is Lipschitz smooth,
    \begin{equation*}
    \Delta_T^{*g}(\mathcal{A}, \mathcal{F}_\mathrm{CSO}^-,\Phi_m^c)\geq \cO( m^{-1} +  \sigma T^{-2/3}).
    \end{equation*}
\end{itemize}
\end{thm}

Together with Theorems \ref{thm:stronglyconvex}, \ref{thm:convex}, and \ref{thm:wc_convergence}, Theorem \ref{thm:lower_bound_CSO} demonstrates that the sample complexity of BSGD cannot be further improved  for strongly convex, convex and weakly convex CSO problems without any additional Lipschitz continuity assumption on the gradient estimator. Similarly, the sample complexity of BSpiderBoost cannot be improved for the nonconvex smooth CSO problems.

\section{Numerical Experiments}
\label{section:numerical}
In this section, we illustrate the performance of the proposed algorithms on the invariant logistic regression and MAML. 
The detailed experiment setup, results, and platform information are deferred to Appendix \ref{secapp:numerical}.

\paragraph{Invariant Logistic Regression}
Invariant learning has wide applications in training robust classifiers \citep{mroueh2015learning, anselmi2013unsupervised}. We consider the invariant logistic regression problem:
\begin{equation}
\label{eq:inv_logistic}
    \min_w \mathbb{E}_{\xi=(a,b)}\big[\log(1+\exp(-b\mathbb{E}_{\eta|\xi}[\eta^Tw])\big],
\end{equation}
where $a\in\mathbb{R}^d$ is the random feature vector, $b\in\{\pm 1\}$ is the corresponding label and $\eta$ is a random perturbed observation of the feature $a$. Let $\sigma_1^2$, $\sigma_2^2$ denote the variances of $a$ and $\eta|a$, respectively. We observe that for a given budget of total samples, BSGD outperforms SAA and converges even when a small inner batch size is used as shown in Table \ref{tab:BSGD_Logistic_Speed}. Detailed results are in Table \ref{tab:BSGD_SAA_Logistic_FULL} in  Appendix \ref{secapp:exp_log}. 

\vspace{-2mm}
\begin{table}[ht]
\caption{Comparison of  BSGD and SAA}
	\renewcommand\arraystretch{1.1}
	\vskip 0.1in
	\small
		\centering
		\begin{tabular}{ccccccc}
			\toprule
			\multirow{2}{*}{$\sigma^2_2/\sigma_1^2$} 
			& \multicolumn{3}{c}{\textbf{BSGD}}  
			& \multicolumn{3}{c}{\textbf{SAA}} \\
			\cline{2-7}
			& \textit{$m$} & \textit{Mean} & \textit{Dev} 
			& \textit{$m$} & \textit{Mean} & \textit{Dev} \\
			\hline
			1
			& 5 & 1.77e-04 & 4.70e-05 & 100& 5.56e-04 & 2.81e-04 \\
			\hline
			10
			& 5 & 3.26e-04 & 1.15e-04 &464 &2.14e-03 & 8.45e-04 \\
			\hline
			100
			& 50 & 1.50e-03 & 6.97e-04 &1000 & 1.12e-02 & 6.42e-04 \\
			\bottomrule
		\end{tabular}
	\label{tab:BSGD_Logistic_Speed}
	\end{table}

\begin{figure*}[t]
\begin{center}
	\includegraphics[width=.32\textwidth,trim=10 10 60 40,clip]{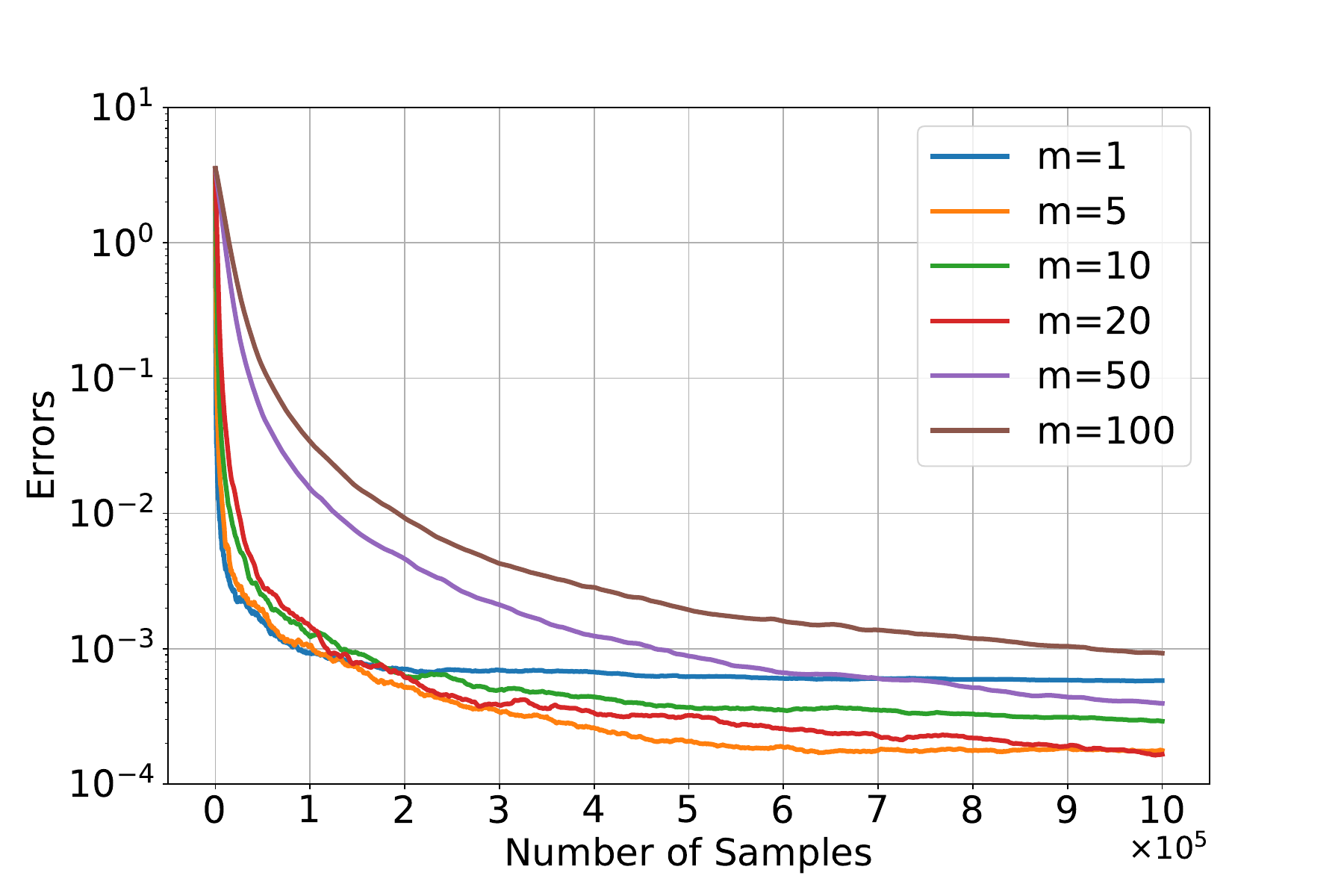}
	\includegraphics[width=.32\textwidth,trim=10 10 60 40,clip]{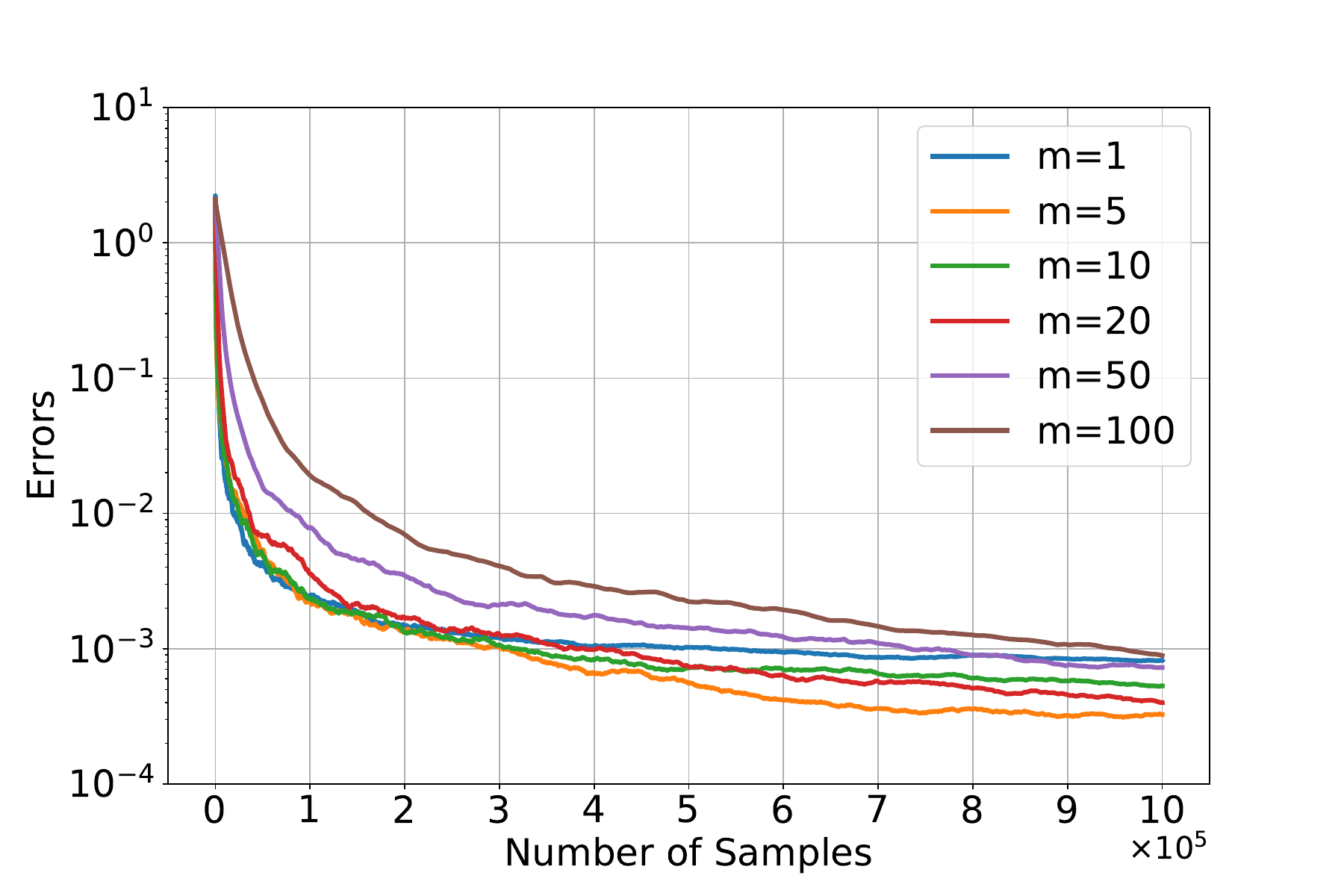}
	\includegraphics[width=.32\textwidth,trim=10 10 60 40,clip]{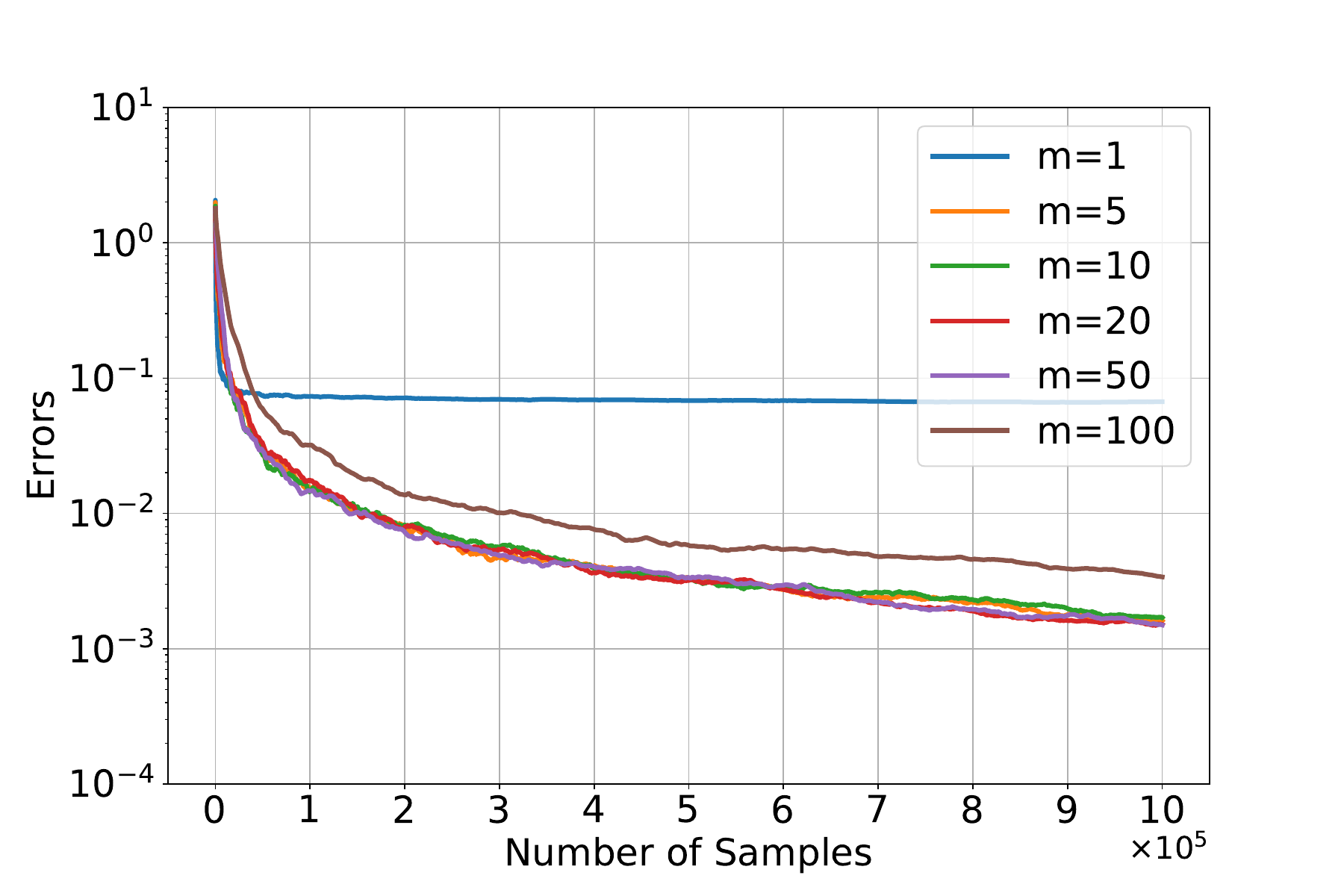}
	\caption{BSGD for invariant Logistic regression (a) $\sigma_2^2/\sigma_1^2=1$, (b) $\sigma_2^2/\sigma_1^2=10$, (c) $\sigma_2^2/\sigma_1^2=100$. }
\label{fig:BSGD-Error-Plot}
\end{center}
 \vskip -0.2in
\end{figure*}

Figure \ref{fig:BSGD-Error-Plot} summarizes the performance of BSGD with different inner batch sizes and under different noise ratios for a given total number of samples. When the noise ratio $\sigma_2^2/\sigma_1^2$ increases, more inner samples are needed to achieve the same performance, as suggested by the theory.

\paragraph{Model-Agnostic Meta-Learning (MAML)}
\label{subsec:app_maml}
We consider the widely used sine-wave few-shot learning. The goal is to find a good initialization model parameter such that the network could recover a new unseen sine wave using only a few available data points. 

The sine wave is of the form $y=a \mathrm{sin}(x+b)$ where $(a,b)$ are drawn from a task distribution. Recall the MAML formulation in \eqref{pro:maml}. In this experiments, we set $\alpha = 0.01$, $l_i(w,D^i) = (y^i-h_i(w,x^i))^2$, where $D^i = (x^i, y^i)$ is the data for the $i$-th task and $h_i$ is a neural network consisting of $2$ hidden layers with $40$ nodes and ReLU activation function between each layers.  We evaluate the MAML objective via empirical objective obtained by empirical risk minimization.  

\begin{figure*}[t]
\centering
\subfigure[]{
\includegraphics[width=.315\textwidth,trim=10 10 80 40,clip]{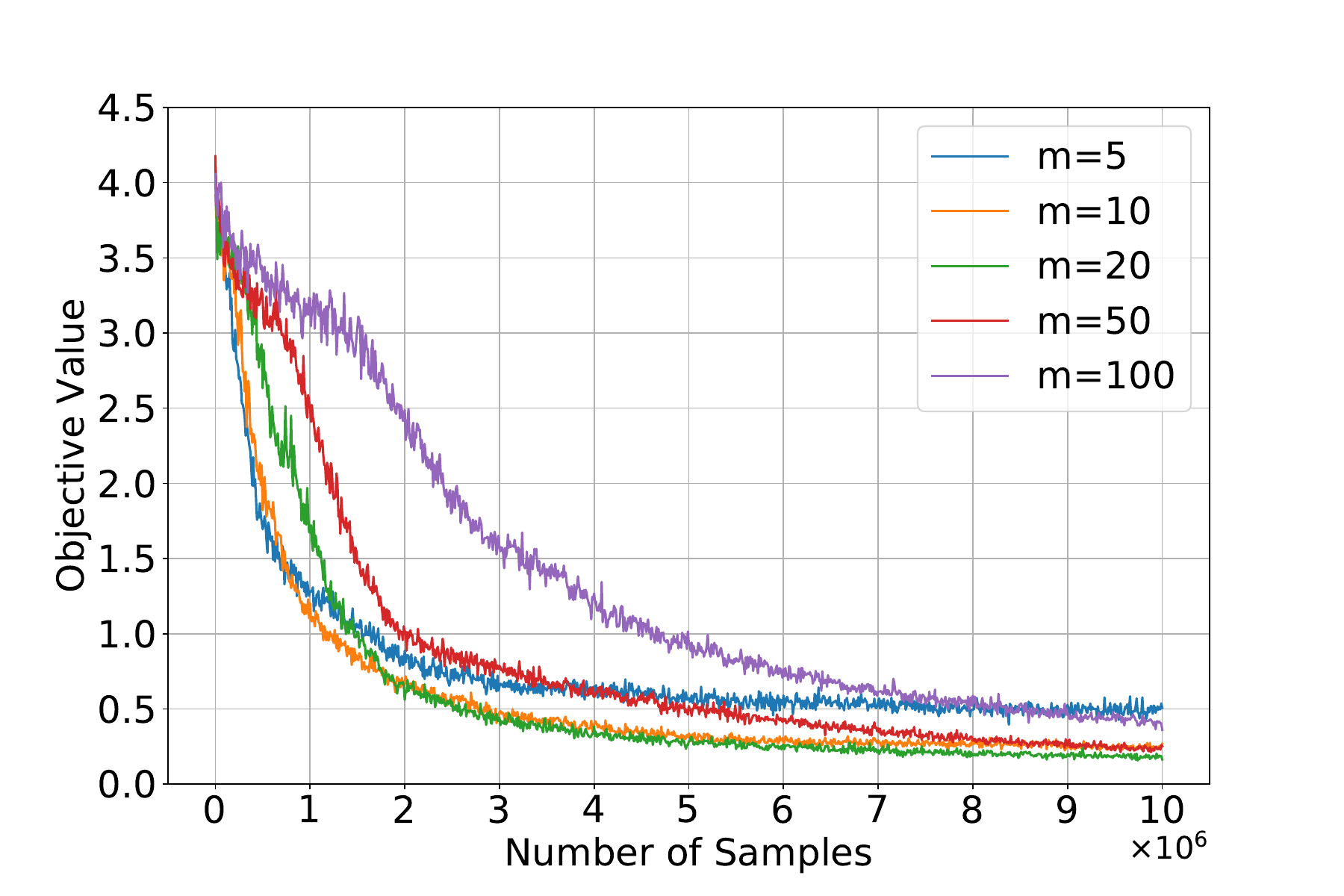}
}
\subfigure[]{
\includegraphics[width=.315\textwidth,trim=10 10 80 40,clip]{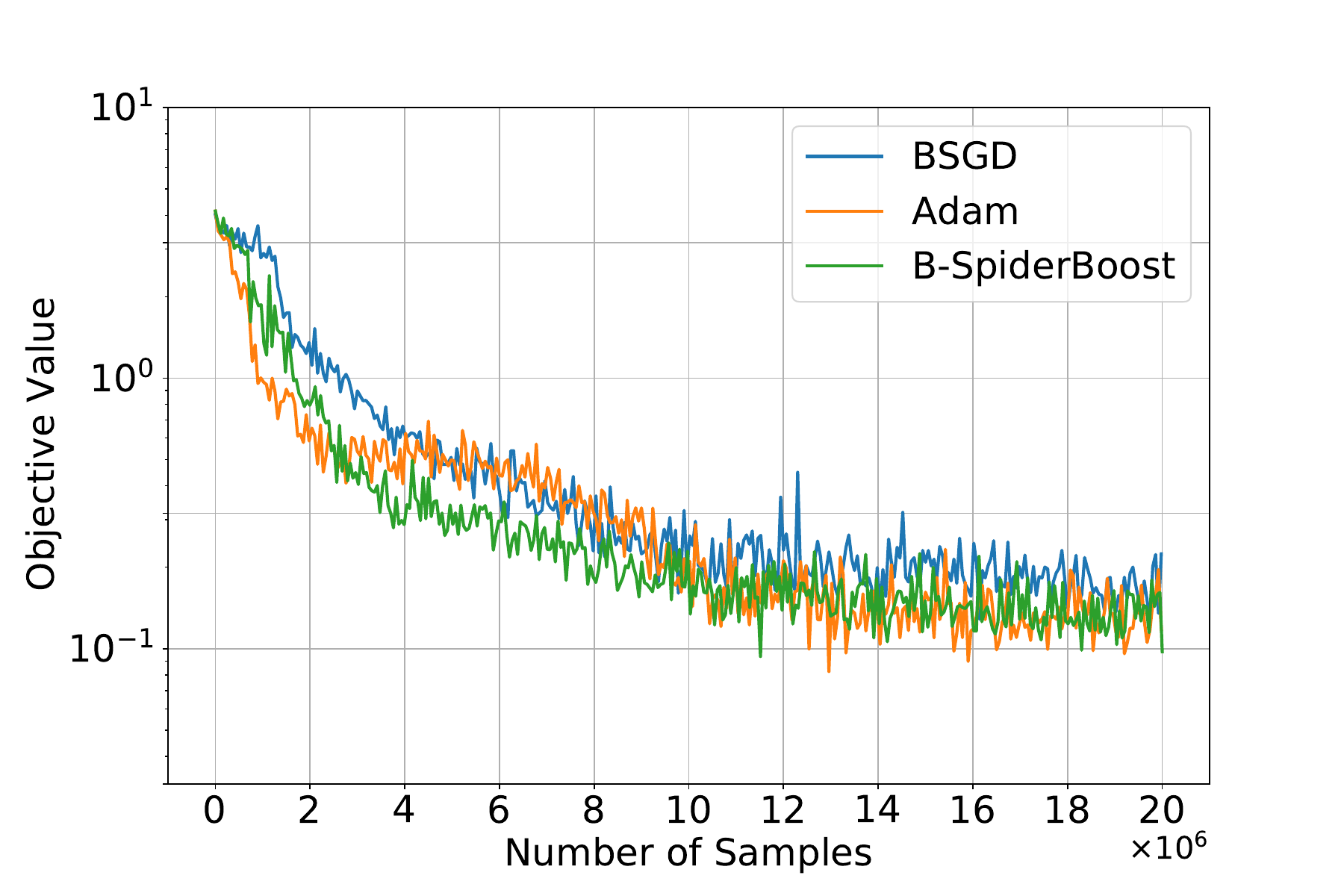}
}
\subfigure[]{
\includegraphics[width=.315\textwidth,trim=10 10 80 40,clip]{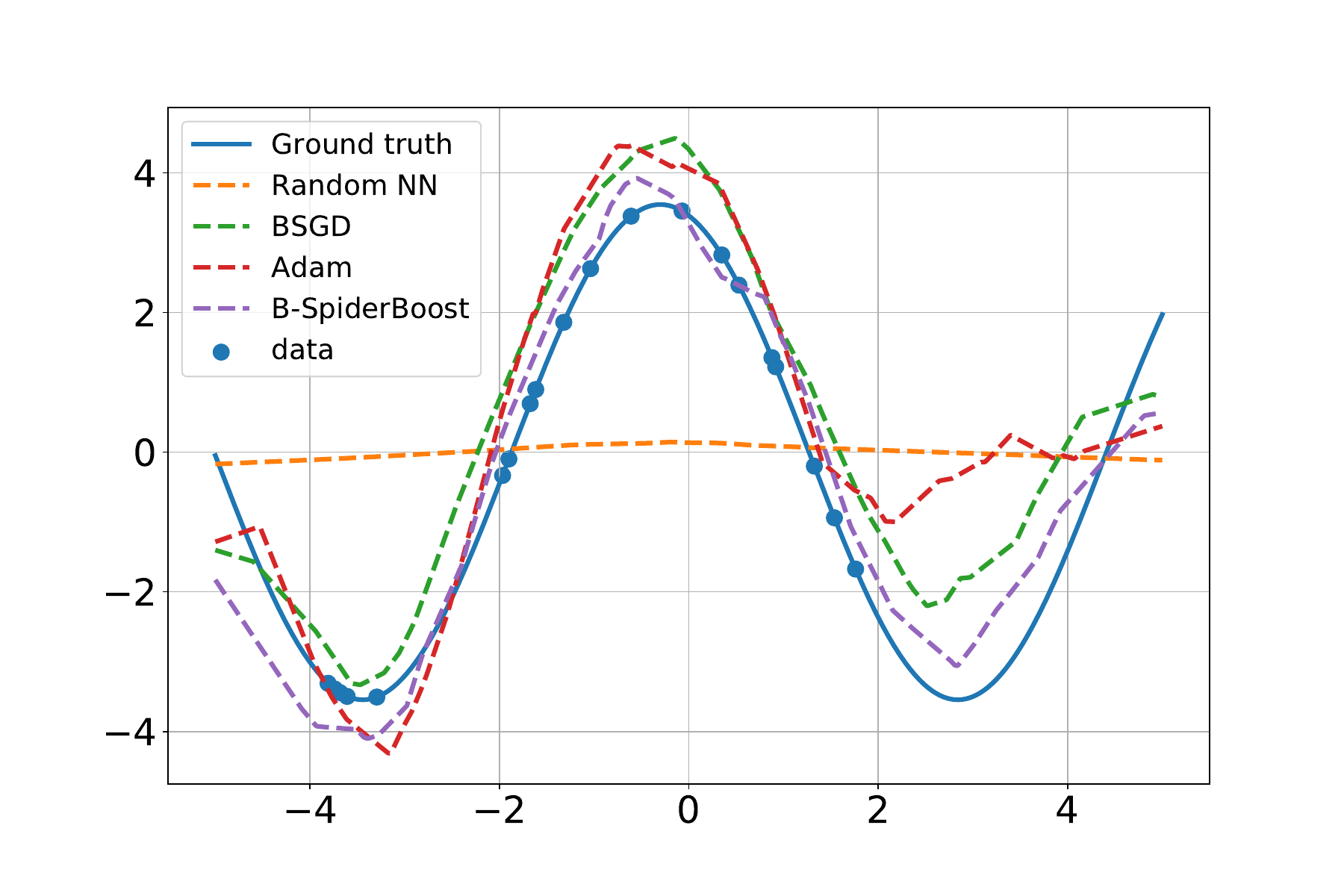}
}
\caption{(a) Convergences of BSGD under differnt inner batch size. (b) Convergences of BSGD, Adam and BSpiderBoost.  (c) Recovered sine-wave signals on an unseen task.}
\label{fig:BSGD-MAML-Plot}
 \vskip -0.2in
\end{figure*}

Figure \ref{fig:BSGD-MAML-Plot}(a)  demonstrates a tradeoff between the inner batch size $m$ and the number of iterations for BSGD. 
Figure \ref{fig:BSGD-MAML-Plot}(b) compares the convergence performance of BSGD, Adam, and BSpiderBoost with the best tuned inner batch sizes. Here Adam refers to a variant of BSGD that performs Adam updates using the biased gradient estimator we constructed. 
Figure \ref{fig:BSGD-MAML-Plot}(c) shows the recovered signal after a one-step update on the unseen task with only $20$ samples using the initialization model parameters obtained by all three algorithms in the meta training step. Random NN refers to the recovered signal using the neural network with random initialization.  

Table \ref{tab:BSGD_MAML_Speed} summarizes the average loss and running time (in CPU minutes) over $10$ trials of each algorithm (under their best inner batch sizes). Although  the widely used first-order MAML (FO-MAML)~\citep{finn2017model} requires the least running time, its performance is worse than BSGD. When $m=50$, FO-MAML does not converge (Figure \ref{fig:FO-MAML} in Appendix).
BSGD requires a smaller batch size to achieve its best performance, which is more practical when a task only has a small number of samples.  
\vspace{-3mm}
\begin{table}[hb]
    \renewcommand\arraystretch{1.1}
    \small
    \caption{Comparison of the average  loss and average running time}
     \vskip 0.1in
    \centering
    \begin{tabular}{ccccccc}
\toprule
        \multicolumn{7}{c}{$\alpha=0.01, Q=10^7$} \\
        \cline{1-7}
        \multirow{2}{*}{\textbf{$m$}} 
        & \multicolumn{2}{c}{\textbf{BSGD}} 
        & \multicolumn{2}{c}{\textbf{FO-MAML}} 
        & \multicolumn{2}{c}{\textbf{Adam}} \\
        \cline{2-7}
        & \textit{Mean} & \textit{CPU} & \textit{Mean} & \textit{CPU} &\textit{Mean} & \textit{CPU} \\
        \hline
        10
        & 2.12e-01 & 71.57 & 2.52e-01 & 41.45 & 8.16e-01 & 86.54 \\
        \hline
        20
        & 2.04e-01 & 35.63 & 2.50e+00 & 20.60 & 3.99e-01 & 43.42 \\
        \hline
        50
        & 2.17e-01 & 14.63 & 3.98e+00 & 8.64 & 2.77e-01 & 17.62 \\
\bottomrule
    \end{tabular}
    \label{tab:BSGD_MAML_Speed}
    \end{table}

\vspace{-3mm}
To summarize, BSpiderBoost achieves the best recovery result but is much harder to tune in practice. In terms of convergence, BSpiderBoost is marginally better than BSGD on this sine-wave task. A possible reason might be that the objective function in our example is not necessarily smooth due to the ReLU activation. Ramdon NN could fail the MAML task when there is a limited amount of samples. More details are available in Appendix \ref{secapp:exp_MAML}.

\clearpage
\section*{Broader Impact}
This paper is purely theoretical and has no immediate ethical or societal consequences.

\begin{ack}
We thank all reviewers and the area chair for the detailed feedback. Funding for this  work was provided by NSF CRII under the award number CCF-1755829.
\end{ack}

\bibliographystyle{plainnat}
\bibliography{main}

\begin{thebibliography}{48}
\providecommand{\natexlab}[1]{#1}
\providecommand{\url}[1]{\texttt{#1}}
\expandafter\ifx\csname urlstyle\endcsname\relax
  \providecommand{\doi}[1]{doi: #1}\else
  \providecommand{\doi}{doi: \begingroup \urlstyle{rm}\Url}\fi

\bibitem[Ajalloeian and Stich(2020)]{ajalloeian2020analysis}
Ahmad Ajalloeian and Sebastian~U Stich.
\newblock Analysis of {SGD} with biased gradient estimators.
\newblock \emph{arXiv preprint arXiv:2008.00051}, 2020.

\bibitem[Anselmi et~al.(2016)Anselmi, Leibo, Rosasco, Mutch, Tacchetti, and Poggio]{anselmi2013unsupervised}
Fabio Anselmi, Joel~Z Leibo, Lorenzo Rosasco, Jim Mutch, Andrea Tacchetti, and Tomaso Poggio.
\newblock Unsupervised learning of invariant representations.
\newblock \emph{Theoretical Computer Science}, 633:\penalty0 112--121, 2016.

\bibitem[Arjevani et~al.(2019)Arjevani, Carmon, Duchi, Foster, Srebro, and Woodworth]{arjevani2019lower}
Yossi Arjevani, Yair Carmon, John~C Duchi, Dylan~J Foster, Nathan Srebro, and Blake Woodworth.
\newblock Lower bounds for non-convex stochastic optimization.
\newblock \emph{arXiv preprint arXiv:1912.02365}, 2019.

\bibitem[Beck(2017)]{beck17firstorder}
Amir Beck.
\newblock \emph{First-Order Methods in Optimization}.
\newblock Society for Industrial and Applied Mathematics, Philadelphia, PA, 2017.

\bibitem[Braun et~al.(2017)Braun, Guzm{\'a}n, and Pokutta]{braun2017lower}
G{\'a}bor Braun, Crist{\'o}bal Guzm{\'a}n, and Sebastian Pokutta.
\newblock Lower bounds on the oracle complexity of nonsmooth convex optimization via information theory.
\newblock \emph{IEEE Transactions on Information Theory}, 63\penalty0 (7):\penalty0 4709--4724, 2017.

\bibitem[Bubeck et~al.(2015)]{bubeck2015convex}
S{\'e}bastien Bubeck et~al.
\newblock Convex optimization: Algorithms and complexity.
\newblock \emph{Foundations and Trends{\textregistered} in Machine Learning}, 8\penalty0 (3-4):\penalty0 231--357, 2015.

\bibitem[Bujok et~al.(2015)Bujok, Hambly, and Reisinger]{bujok2015multilevel}
Karolina Bujok, BM~Hambly, and Christoph Reisinger.
\newblock Multilevel simulation of functionals of bernoulli random variables with application to basket credit derivatives.
\newblock \emph{Methodology and Computing in Applied Probability}, 17\penalty0 (3):\penalty0 579--604, 2015.

\bibitem[Chen and Luss(2018)]{chen2018stochastic}
Jie Chen and Ronny Luss.
\newblock Stochastic gradient descent with biased but consistent gradient estimators.
\newblock \emph{arXiv preprint arXiv:1807.11880}, 2018.

\bibitem[Chen et~al.(2020)Chen, Sun, and Yin]{chen2020solving}
Tianyi Chen, Yuejiao Sun, and Wotao Yin.
\newblock Solving stochastic compositional optimization is nearly as easy as solving stochastic optimization.
\newblock \emph{arXiv preprint arXiv:2008.10847}, 2020.

\bibitem[Chen et~al.(2019)Chen, Yuan, Yi, Zhou, Chen, and Yang]{chen2018universal}
Zaiyi Chen, Zhuoning Yuan, Jinfeng Yi, Bowen Zhou, Enhong Chen, and Tianbao Yang.
\newblock Universal stagewise learning for non-convex problems with convergence on averaged solutions.
\newblock In \emph{International Conference on Learning Representations}, 2019.

\bibitem[Cover(1999)]{cover1999elements}
Thomas~M Cover.
\newblock \emph{Elements of information theory}.
\newblock John Wiley \& Sons, 1999.

\bibitem[Dai et~al.(2017)Dai, He, Pan, Boots, and Song]{pmlr-v54-dai17a}
Bo~Dai, Niao He, Yunpeng Pan, Byron Boots, and Le~Song.
\newblock Learning from conditional distributions via dual embeddings.
\newblock In \emph{Artificial Intelligence and Statistics}, pages 1458--1467, 2017.

\bibitem[Dai et~al.(2018)Dai, Shaw, Li, Xiao, He, Liu, Chen, and Song]{pmlr-v80-dai18c}
Bo~Dai, Albert Shaw, Lihong Li, Lin Xiao, Niao He, Zhen Liu, Jianshu Chen, and Le~Song.
\newblock {SBEED}: Convergent reinforcement learning with nonlinear function approximation.
\newblock In \emph{International Conference on Machine Learning}, pages 1125--1134, 2018.

\bibitem[Davis and Drusvyatskiy(2019)]{davis2019stochastic}
Damek Davis and Dmitriy Drusvyatskiy.
\newblock Stochastic model-based minimization of weakly convex functions.
\newblock \emph{SIAM Journal on Optimization}, 29\penalty0 (1):\penalty0 207--239, 2019.

\bibitem[Davis et~al.(2020)Davis, Drusvyatskiy, Kakade, and Lee]{davis2020stochastic}
Damek Davis, Dmitriy Drusvyatskiy, Sham Kakade, and Jason~D Lee.
\newblock Stochastic subgradient method converges on tame functions.
\newblock \emph{Foundations of computational mathematics}, 20\penalty0 (1):\penalty0 119--154, 2020.

\bibitem[Drusvyatskiy(2017)]{drusvyatskiy2017proximal}
Dmitriy Drusvyatskiy.
\newblock The proximal point method revisited.
\newblock \emph{arXiv preprint arXiv:1712.06038}, 2017.

\bibitem[Fallah et~al.(2020)Fallah, Mokhtari, and Ozdaglar]{fallah2019convergence}
Alireza Fallah, Aryan Mokhtari, and Asuman Ozdaglar.
\newblock On the convergence theory of gradient-based model-agnostic meta-learning algorithms.
\newblock In \emph{International Conference on Artificial Intelligence and Statistics}, pages 1082--1092, 2020.

\bibitem[Fang et~al.(2018)Fang, Li, Lin, and Zhang]{fang2018spider}
Cong Fang, Chris~Junchi Li, Zhouchen Lin, and Tong Zhang.
\newblock Spider: Near-optimal non-convex optimization via stochastic path-integrated differential estimator.
\newblock In \emph{Advances in Neural Information Processing Systems}, pages 689--699, 2018.

\bibitem[Finn et~al.(2017)Finn, Abbeel, and Levine]{finn2017model}
Chelsea Finn, Pieter Abbeel, and Sergey Levine.
\newblock Model-agnostic meta-learning for fast adaptation of deep networks.
\newblock In \emph{Proceedings of the 34th International Conference on Machine Learning-Volume 70}, pages 1126--1135. JMLR. org, 2017.

\bibitem[Ghadimi et~al.(2020)Ghadimi, Ruszczynski, and Wang]{ghadimi2018single}
Saeed Ghadimi, Andrzej Ruszczynski, and Mengdi Wang.
\newblock A single timescale stochastic approximation method for nested stochastic optimization.
\newblock \emph{SIAM Journal on Optimization}, 30\penalty0 (1):\penalty0 960--979, 2020.

\bibitem[Giles(2008)]{giles2008multilevel}
Michael~B Giles.
\newblock Multilevel monte carlo path simulation.
\newblock \emph{Operations research}, 56\penalty0 (3):\penalty0 607--617, 2008.

\bibitem[Giles(2015)]{giles2015multilevel}
Michael~B Giles.
\newblock Multilevel monte carlo methods.
\newblock \emph{Acta Numerica}, 24:\penalty0 259, 2015.

\bibitem[Giles and Haji-Ali(2019)]{giles2019multilevel}
Michael~B Giles and Abdul-Lateef Haji-Ali.
\newblock Multilevel nested simulation for efficient risk estimation.
\newblock \emph{SIAM/ASA Journal on Uncertainty Quantification}, 7\penalty0 (2):\penalty0 497--525, 2019.

\bibitem[Gordy and Juneja(2010)]{gordy2010nested}
Michael~B Gordy and Sandeep Juneja.
\newblock Nested simulation in portfolio risk measurement.
\newblock \emph{Management Science}, 56\penalty0 (10):\penalty0 1833--1848, 2010.

\bibitem[Hazan et~al.(2007)Hazan, Agarwal, and Kale]{hazan2007logarithmic}
Elad Hazan, Amit Agarwal, and Satyen Kale.
\newblock Logarithmic regret algorithms for online convex optimization.
\newblock \emph{Machine Learning}, 69\penalty0 (2-3):\penalty0 169--192, 2007.

\bibitem[Hong and Juneja(2009)]{hong2009estimating}
L~Jeff Hong and Sandeep Juneja.
\newblock Estimating the mean of a non-linear function of conditional expectation.
\newblock In \emph{Winter Simulation Conference}, pages 1223--1236. Winter Simulation Conference, 2009.

\bibitem[Hong et~al.(2017)Hong, Juneja, and Liu]{hong2017kernel}
L~Jeff Hong, Sandeep Juneja, and Guangwu Liu.
\newblock Kernel smoothing for nested estimation with application to portfolio risk measurement.
\newblock \emph{Operations Research}, 65\penalty0 (3):\penalty0 657--673, 2017.

\bibitem[Hu et~al.(2020{\natexlab{a}})Hu, Seiler, and Lessard]{hu2020analysis}
Bin Hu, Peter Seiler, and Laurent Lessard.
\newblock Analysis of biased stochastic gradient descent using sequential semidefinite programs.
\newblock \emph{Mathematical Programming}, pages 1--26, 2020{\natexlab{a}}.

\bibitem[Hu et~al.(2016)Hu, Prashanth, Gy{\"o}rgy, and Szepesv{\'a}ri]{hu2016bandit}
Xiaowei Hu, LA~Prashanth, Andr{\'a}s Gy{\"o}rgy, and Csaba Szepesv{\'a}ri.
\newblock (bandit) convex optimization with biased noisy gradient oracles.
\newblock In \emph{Artificial Intelligence and Statistics}, pages 819--828, 2016.

\bibitem[Hu et~al.(2020{\natexlab{b}})Hu, Chen, and He]{hu2019sample}
Yifan Hu, Xin Chen, and Niao He.
\newblock Sample complexity of sample average approximation for conditional stochastic optimization.
\newblock \emph{SIAM Journal on Optimization}, 30\penalty0 (3):\penalty0 2103--2133, 2020{\natexlab{b}}.

\bibitem[Johnson and Zhang(2013)]{johnson2013accelerating}
Rie Johnson and Tong Zhang.
\newblock Accelerating stochastic gradient descent using predictive variance reduction.
\newblock In \emph{Advances in neural information processing systems}, pages 315--323, 2013.

\bibitem[Karimi et~al.(2019)Karimi, Miasojedow, Moulines, and Wai]{karimi2019non}
Belhal Karimi, Blazej Miasojedow, Eric Moulines, and Hoi-To Wai.
\newblock Non-asymptotic analysis of biased stochastic approximation scheme.
\newblock In \emph{Conference on Learning Theory}, pages 1944--1974, 2019.

\bibitem[Mroueh et~al.(2015)Mroueh, Voinea, and Poggio]{mroueh2015learning}
Youssef Mroueh, Stephen Voinea, and Tomaso~A Poggio.
\newblock Learning with group invariant features: A kernel perspective.
\newblock In \emph{Advances in Neural Information Processing Systems}, pages 1558--1566, 2015.

\bibitem[Muandet et~al.(2019)Muandet, Mehrjou, Lee, and Raj]{muandet2019dual}
Krikamol Muandet, Arash Mehrjou, Si~Kai Lee, and Anant Raj.
\newblock Dual {IV}: A single stage instrumental variable regression.
\newblock \emph{arXiv preprint arXiv:1910.12358}, 2019.

\bibitem[Nachum and Dai(2020)]{nachum2020}
Ofir Nachum and Bo~Dai.
\newblock Reinforcement learning via fenchel-rockafellar duality.
\newblock \emph{arXiv preprint arXiv:2001.01866}, 2020.

\bibitem[Nemirovski et~al.(2009)Nemirovski, Juditsky, Lan, and Shapiro]{nemirovski2009robust}
Arkadi Nemirovski, Anatoli Juditsky, Guanghui Lan, and Alexander Shapiro.
\newblock Robust stochastic approximation approach to stochastic programming.
\newblock \emph{SIAM Journal on optimization}, 19\penalty0 (4):\penalty0 1574--1609, 2009.

\bibitem[Nguyen et~al.(2017)Nguyen, Liu, Scheinberg, and Tak{\'a}{\v{c}}]{nguyen2017sarah}
Lam~M Nguyen, Jie Liu, Katya Scheinberg, and Martin Tak{\'a}{\v{c}}.
\newblock Sarah: A novel method for machine learning problems using stochastic recursive gradient.
\newblock In \emph{Proceedings of the 34th International Conference on Machine Learning-Volume 70}, pages 2613--2621. JMLR. org, 2017.

\bibitem[Pang et~al.(2017)Pang, Razaviyayn, and Alvarado]{pang2017computing}
Jong-Shi Pang, Meisam Razaviyayn, and Alberth Alvarado.
\newblock Computing b-stationary points of nonsmooth dc programs.
\newblock \emph{Mathematics of Operations Research}, 42\penalty0 (1):\penalty0 95--118, 2017.

\bibitem[Paquette et~al.(2018)Paquette, Lin, Drusvyatskiy, Mairal, and Harchaoui]{paquette2017catalyst}
Courtney Paquette, Hongzhou Lin, Dmitriy Drusvyatskiy, Julien Mairal, and Zaid Harchaoui.
\newblock Catalyst for gradient-based nonconvex optimization.
\newblock In \emph{International Conference on Artificial Intelligence and Statistics}, volume~84 of \emph{Proceedings of Machine Learning Research}, pages 613--622. PMLR, 2018.

\bibitem[Shamir and Zhang(2013)]{shamir2013stochastic}
Ohad Shamir and Tong Zhang.
\newblock Stochastic gradient descent for non-smooth optimization: Convergence results and optimal averaging schemes.
\newblock In \emph{International Conference on Machine Learning}, pages 71--79, 2013.

\bibitem[Singh et~al.(2019)Singh, Sahani, and Gretton]{singh2019kernel}
Rahul Singh, Maneesh Sahani, and Arthur Gretton.
\newblock Kernel instrumental variable regression.
\newblock In \emph{Advances in Neural Information Processing Systems 32}, pages 4595--4607. Curran Associates, Inc., 2019.

\bibitem[Wang et~al.(2016)Wang, Liu, and Fang]{wang2016accelerating}
Mengdi Wang, Ji~Liu, and Ethan Fang.
\newblock Accelerating stochastic composition optimization.
\newblock In \emph{Advances in Neural Information Processing Systems}, pages 1714--1722, 2016.

\bibitem[Wang et~al.(2017)Wang, Fang, and Liu]{wang2017stochastic}
Mengdi Wang, Ethan~X Fang, and Han Liu.
\newblock Stochastic compositional gradient descent: algorithms for minimizing compositions of expected-value functions.
\newblock \emph{Mathematical Programming}, 161\penalty0 (1-2):\penalty0 419--449, 2017.

\bibitem[Wang et~al.(2019)Wang, Ji, Zhou, Liang, and Tarokh]{wang2018spiderboost}
Zhe Wang, Kaiyi Ji, Yi~Zhou, Yingbin Liang, and Vahid Tarokh.
\newblock Spiderboost and momentum: Faster variance reduction algorithms.
\newblock In \emph{Advances in Neural Information Processing Systems}, pages 2406--2416, 2019.

\bibitem[Yang et~al.(2019)Yang, Wang, and Fang]{yang2019multilevel}
Shuoguang Yang, Mengdi Wang, and Ethan~X Fang.
\newblock Multilevel stochastic gradient methods for nested composition optimization.
\newblock \emph{SIAM Journal on Optimization}, 29\penalty0 (1):\penalty0 616--659, 2019.

\bibitem[{Yao}(1977)]{yao1977}
A.~C. {Yao}.
\newblock Probabilistic computations: Toward a unified measure of complexity.
\newblock In \emph{18th Annual Symposium on Foundations of Computer Science (sfcs 1977)}, pages 222--227, Oct 1977.

\bibitem[Zhang and Xiao(2019)]{zhang2019multi}
Junyu Zhang and Lin Xiao.
\newblock Multi-level composite stochastic optimization via nested variance reduction.
\newblock \emph{arXiv preprint arXiv:1908.11468}, 2019.

\bibitem[Zhang and He(2018)]{zhang2018convergence}
Siqi Zhang and Niao He.
\newblock On the convergence rate of stochastic mirror descent for nonsmooth nonconvex optimization.
\newblock \emph{arXiv preprint arXiv:1806.04781}, 2018.

\end{thebibliography}

\newpage
\appendix
\onecolumn
\begin{appendix}
\begin{center}
{\huge Appendix}
\end{center}
\section{Proof of Lemma \ref{lm:bias_gradient}}
\begin{proof}
Denote $\hat g(x,\xi) := \frac{1}{m}\sum_{j=1}^m  g_{\eta_j}(x,\xi)$.
Note that
\begin{equation*}
	\begin{split}
	&\|\EE \nabla \hat{F}(x; \xi,\{\eta_{j}\}_{j=1}^{m}) 
	-
	\nabla F(x)\|_2^2\\
	\leq\ &
	\left\|\EE_\xi \EE_{\{\eta_j\}_{j=1}^m|\xi} \left(\nabla \hat g(x,\xi)\right)^\top \nabla f_\xi\left(\hat g(x,\xi)\right)
	-
	\EE_\xi (\EE_{\eta|\xi}\nabla g_\eta(x,\xi))^\top \nabla f_\xi(\EE_{\eta|\xi}g_\eta(x,\xi) ) \right\|_2^2\\
	\leq \ & 
	\left\|\EE_\xi \EE_{\{\eta_j\}_{j=1}^m|\xi} \left(\nabla \hat g(x,\xi)\right)^\top \nabla f_\xi\left(\hat g(x,\xi)\right)
	-
	\EE_\xi \EE_{\{\eta_j\}_{j=1}^m|\xi} \left(\nabla \hat g(x,\xi)\right)^\top \nabla f_\xi\left(\EE_{\eta|\xi} g_\eta(x,\xi)\right)\right\|_2^2\\
	+\ & \left\|\EE_\xi \EE_{\{\eta_j\}_{j=1}^m|\xi} \left(\nabla \hat g(x,\xi)\right)^\top \nabla f_\xi\left(\EE_{\eta|\xi} g_\eta(x,\xi)\right)
	-
	\EE_\xi (\EE_{\eta|\xi}\nabla g_\eta(x,\xi))^\top \nabla f_\xi(\EE_{\eta|\xi}g_\eta(x,\xi) ) \right\|_2^2\\
	= \ & 
	\left\|\EE_\xi \EE_{\{\eta_j\}_{j=1}^m|\xi} \left(\nabla \hat g(x,\xi)\right)^\top \nabla f_\xi\left(\hat g(x,\xi)\right)
	-
	\EE_\xi \EE_{\{\eta_j\}_{j=1}^m|\xi} \left(\nabla \hat g(x,\xi)\right)^\top \nabla f_\xi\left(\EE_{\eta|\xi} g_\eta(x,\xi)\right)\right\|_2^2\\
	\leq\ &
	\EE_\xi \EE_{\{\eta_j\}_{j=1}^m|\xi} \|\nabla \hat g(x,\xi)\|_2^2 \|\nabla f_\xi\left(\hat g(x,\xi)\right)
	-  
	\nabla f_\xi\left(\EE_{\eta|\xi} g_\eta(x,\xi)\right)\|_2^2\\
	\leq \ & 
	L_g^2 S_f^2\EE_\xi \EE_{\{\eta_j\}_{j=1}^m|\xi} \|\hat g(x,\xi)-\EE_{\eta|\xi} g_\eta(x,\xi) \|_2^2\\
	\leq \ & 
	\frac{L_g^2S_f^2\sigma_g^2}{m}.
	\end{split}
	\end{equation*}
The equality holds as 
\begin{equation*}
\begin{split}
& \EE_\xi \EE_{\{\eta_j\}_{j=1}^m|\xi} \left(\nabla \hat g(x,\xi)\right)^\top \nabla f_\xi\left(\EE_{\eta|\xi} g_\eta(x,\xi)\right)
-
\EE_\xi (\EE_{\eta|\xi}\nabla g_\eta(x,\xi))^\top \nabla f_\xi(\EE_{\eta|\xi}g_\eta(x,\xi) ) \\
= & \EE_\xi \EE_{\{\eta_j\}_{j=1}^m|\xi} (\nabla \hat g(x,\xi)- \EE_{\eta|\xi}\nabla g_\eta(x,\xi))^\top \nabla f_\xi(\EE_{\eta|\xi}g_\eta(x,\xi) )\\
=  & 0.   
\end{split}
\end{equation*}
\end{proof}

\section{Convergence Analysis}
\label{secapp:convergence}
In this section, we present the proof of Theorems \ref{thm:stronglyconvex}, \ref{thm:convex}, and \ref{thm:wc_convergence}. Based on these theorems, we demonstrate the sample complexity of BSGD with strongly convex, convex, and weakly convex objectives.

First, we present the proof framework for strongly convex and convex objectives. Recall  BSGD in Algorithm \ref{alg:A} , at iteration $t$, BSGD first generates sample $\xi_t$  from the distribution of $\xi$ and $m$ samples $\{\eta_{tj}\}_{j=1}^{m_t}$ from the conditional distribution of $\eta|\xi_t$. We define the following auxiliary functions to facilitate our analysis:
\begin{equation*}
p(x,\xi_t) := f_{\xi_t}(\EE_{\eta|\xi_t} g_\eta(x,\xi_t)); \quad \hat p(x,\xi_t): =f_{\xi_t}\Big(\frac{1}{m_t}\sum_{j=1}^{m_t} g_{\eta_{tj}}(x,{\xi_t})\Big).    
\end{equation*}
Note that $\hat F(x;\xi_t,\{\eta_{tj}\}_{j=1}^{m_t})=\hat p(x,\xi_t)$. The biased gradient estimator used in BSGD is $\nabla \hat p(x,\xi_t)$.
Denote $x^*\in \argmin_{x\in\mathcal{X}} F(x)$, $A_t = \frac{1}{2}\norm{x_t-x^*}_2^2$, $a_t = \EE A_t$.
Since $\Pi_\mathcal{X}(x^*)=x^*$ and the projection operator is non-expansive, we have
\begin{equation}
\begin{split}
    A_{t+1} 
    & =  \frac{1}{2}\norm{x_{t+1}-x^*}_2^2 \\
    & =  \frac{1}{2}\norm{\Pi_\mathcal{X}(x_{t}-\gamma_t \nabla_{x} \hat p(x_t,\xi_t))-\Pi_\mathcal{X}(x^*)}_2^2 \\
    & \leq \frac{1}{2}\norm{x_{t}-x^*-\gamma_t \nabla_{x} \hat p(x_t,\xi_t)}_2^2\\
    & = A_t + \frac{1}{2}\gamma_t^2\norm{\nabla_{x} \hat p(x_t,\xi_t)}_2^2 - \gamma_t \nabla_{x} \hat p(x_t,\xi_t)^\top(x_{t}-x^*).
\end{split}
\end{equation}
Dividing $\gamma_t$ on both sides and taking expectation over $\{\xi_t,\{\eta_{tj}\}_{j=1}^{m_t}\}$, it holds
\begin{equation}
\label{eq:point}
    \EE \nabla_{x} \hat p(x_t,\xi_t)^\top(x_{t}-x^*)\leq \frac{a_t-a_{t+1}}{\gamma_t}  + \frac{1}{2}\gamma_t\EE \norm{\nabla_{x} \hat p(x_t,\xi_t)}_2^2.
\end{equation}
By Assumption \ref{ass:convex},  we have
\begin{equation}
\begin{split}
\label{eq:mu_convex_on_p}
    &-\nabla_{x} \hat p(x_t,\xi_t)^\top (x_t-x^*) 
    \leq \;  \hat p(x^*,\xi_t)-\hat p(x_t,\xi_t)-\frac{\mu}{2}\|x_t-x^*\|_2^2\\
    = \; & 
    \underbrace{\hat p(x^*,\xi_t)- p(x^*,\xi_t)}_{:=\zeta_{t1}} 
    + \underbrace{p(x^*,\xi_t)- p(x_t,\xi_t)}_{:=\zeta_{t2}}
    + \underbrace{p(x_t,\xi_t) - \hat p(x_t,\xi_t)}_{:=\zeta_{t3}}-\frac{\mu}{2}\|x_t-x^*\|_2^2.
\end{split}
\end{equation}
Taking expectation over $\{\xi_t,\{\eta_{tj}\}_{j=1}^{m_t}\}$ on both sides, by the definition of $p(x,\xi)$, it holds $\EE_{\xi_t}[ \zeta_{t2}\mid x_t]=\EE_{\xi_t} [ F(x^*)-F(x_t)\mid x_t]$, then 
\begin{equation}
\label{eq:convex}
-\EE\nabla_{x} \hat p(x_t,\xi_t)^\top (x_t-x^*)\leq \EE \zeta_{t1}+\EE\zeta_{t3}+ \EE [F(x^*)- F(x_t)]-\mu a_t.
\end{equation}
Since $x^*$ and $x_t$ are independent of $\{\xi_t,\{\eta_{tj}\}_{j=1}^m\}$, by Lemma \ref{lm:bias}, we upper bound  $\EE \zeta_{t1}$ and $\EE\zeta_{t3}$ using $\Delta_f(m_t)$:
\begin{equation*}
|\EE \zeta_{t1}|\leq \Delta_f(m_t), \quad |\EE \zeta_{t3}|\leq \Delta_f(m_t).
\end{equation*}
Summing up \eqref{eq:point}  and \eqref{eq:convex}, we obtain
\begin{equation}
\label{eq:key}
\EE [F(x_t)-F(x^*)]
\leq
2\Delta_f(m_t)-\mu a_t+\frac{a_t-a_{t+1}}{\gamma_t}  + \frac{1}{2}\gamma_t\EE \norm{\nabla_{x} \hat p(x_t,\xi_t)}_2^2.
\end{equation}
By convexity of $F$ and the definition of $\hat{x}_T = \frac{1}{T}\sum_{t=1}^T x_t$, we have
\begin{equation}
\label{eq:key2}
\EE[F(\hat{x}_T)-F(x^*)] 
= 
\EE \bigg[F\bigg(\frac{1}{T}\sum_{t=1}^T x_t\bigg)-F(x^*)\bigg]
\leq 
\frac{1}{T} \sum_{t=1}^T\EE[F(x_t)-F(x^*)],    
\end{equation}
We then prove the convergence of BSGD for strongly convex and convex objectives based on \eqref{eq:key}.

\subsection{Global Convergence of BSGD for Strongly Convex Objectives} 
We prove Theorem \ref{thm:stronglyconvex}, the strongly convex case for which Assumption \ref{ass:convex} holds with $\mu>0$. 

\begin{proof}
Since $\hat F$ is $S_F$-Lipschitz smooth and $\mu$-strongly convex, we have
\begin{equation}
\begin{split}
\EE \|\nabla \hat F(x)\|_2^2 \leq & 2\EE\|\nabla \hat F(x)-\nabla \hat F(x^*)\|_2^2+2\EE\|\nabla \hat F(x^*)\|_2^2\\
\leq  & 2S^2 \|x-x^*\|_2^2 +2\EE\|\nabla \hat F(x^*)\|_2^2\\
\leq & 4S^2/\mu (F(x)-F(x^*))+2\EE\|\nabla \hat F(x^*)\|_2^2. \\
\end{split}
\end{equation}
It implies that 
\begin{equation}
\EE [F(x_t)-F(x^*)]
\leq
2\Delta_f(m_t)-\mu a_t+\frac{a_t-a_{t+1}}{\gamma_t}  + \frac{1}{2}\gamma_t(4S_F^2/\mu (F(x)-F(x^*))+2\EE\|\nabla \hat F(x^*)\|_2^2).
\end{equation}
Therefore, we have for $\gamma_t\leq \frac{\mu}{4S_F^2}$, 
\begin{equation}
\label{eq:sc_key}
\begin{split}
\EE F(x_t)-F(x^*) 
& \leq \frac{1}{1-\gamma_t S_F^2/\mu}\bigg( 2\Delta_f(m_t)-\mu a_t+\frac{a_t-a_{t+1}}{\gamma_t}  + \gamma_t\EE\|\nabla \hat F(x^*)\|_2^2\bigg)   \\
& \leq 2\bigg( 2\Delta_f(m_t)-\mu a_t+\frac{a_t-a_{t+1}}{\gamma_t}  + \gamma_t\EE\|\nabla \hat F(x^*)\|_2^2\bigg).
\end{split}
\end{equation}
Summing up \eqref{eq:sc_key} from $t=1$ to $T$,
\begin{equation}
\begin{split}
&\frac{1}{T} \sum_{t=1}^T\EE[F(x_t)-F(x^*)]\\
\leq &
\frac{2}{T}\sum_{t=1}^T\left[ 2\Delta_f(m_t)-\mu a_t+\frac{a_t-a_{t+1}}{\gamma_t}  + \gamma_t\EE\|\nabla \hat F(x^*)\|_2^2\right]\\
\leq & 
\frac{2}{T}\sum_{t=1}^T\left[ 2\Delta_f(m_t)+\gamma_t\EE\|\nabla \hat F(x^*)\|_2^2\right] + \frac{2}{T}\sum_{t=2}^T a_t\left(\frac{1}{\gamma_{t}} -\frac{1}{\gamma_{t-1}}-\mu\right) +\frac{2}{T}a_1\left(\frac{1}{\gamma_{1}} -\mu\right).
\end{split}
\end{equation}
Set $\gamma_t  = \frac{1}{\mu(t+c)}$ and $c =\max\{ 4S_F^2/\mu^2-1,0\}$. It makes sure that $\gamma_t\leq \gamma_1 \leq \frac{\mu}{4S_F^2}$.
Since $1/\gamma_1 -\mu \leq \mu(4S_F^2/\mu^2-1)$, with inequality \eqref{eq:key2}, it holds
\begin{equation*}
\EE[F(\hat{x}_T)-F(x^*)]\leq \frac{4}{T}\sum_{t=1}^T\Delta_f(m_t)+\frac{1}{T}\sum_{t=1}^T \frac{2\EE\|\nabla \hat F(x^*)\|_2^2}{\mu (t+c)} +\frac{ S_F^2}{4\mu T}\|x_1-x^*\|_2^2.
\end{equation*}
By the fact that $\sum_{t=1}^T \frac{1}{t+c} \leq \sum_{t=1}^T  \frac{1}{t}\leq\log(T)+1$, it holds
\begin{equation*}
\EE[F(\hat{x}_T)-F(x^*)]\leq \frac{4}{T}\sum_{t=1}^T\Delta_f(m_t)+\frac{2\EE\|\nabla \hat F(x^*)\|_2^2(\log(T)+1)+ S_F^2/4\|x_1-x^*\|_2^2}{T\mu}.
\end{equation*}
\end{proof} 


We demonstrate the sample complexity using the following corollary.

\begin{cor}
\label{cor:strongly_convex}
To achieve an $\eps$-optimal solution, the total sample complexity of BSGD in the strongly convex case is $\tilde \cO(\eps^{-3})$ for objectives with Lipschitz continuous $f_\xi$ and $\tilde \cO(\eps^{-2})$ for objectives with Lipschitz smooth $f_\xi$. 
\end{cor}
It implies that the smoothness of the outer function makes a difference in the total sample complexity of BSGD when solving CSO.  It is worth pointing out that the sample complexity of BSGD matches with that of ERM (SAA) for strongly convex objectives established in~\citet{hu2019sample}. We now prove Corollary \ref{cor:strongly_convex}.

\begin{proof}
For fixed mini-batch size, to guarantee that $\EE[F(\hat{x}_T)-F(x^*)]\leq \eps$, setting $T=\tilde\cO(\eps^{-1})$ and picking $m = \cO(\eps^{-2})$ for objectives with Lipschitz continuous outer function $f_\xi$ and $m=\cO(\eps^{-1})$ for objectives with Lipschitz smooth outer function $f_\xi$ are sufficient to guarantee that $\hat{x}_T$ is an $\eps$-optimal solution to the \eqref{pro:ori}.

As for time-varying mini-batch sizes, letting $m_t = t^2$ for Lipschitz continuous $f_\xi$. Since $\sum_{t=1}^T \frac{1}{t}\leq \log(T)+1$, it holds
\begin{equation*}
\frac{1}{T}\sum_{t=1}^T\Delta_f(m_t) = \frac{1}{T}\sum_{t=1}^T \frac{L_f\sigma_g}{t}\leq \frac{L_f\sigma_g (\log(T)+1)}{T}\leq \cO(\eps).
\end{equation*}
As a result, setting $T=\tilde \cO(\eps^{-1})$; the total sample complexity is $\sum_{t=1}^T (m_t+1) =\cO(T^3) =\tilde\cO(\eps^{-3})$. 

Set $m_t = t$ for Lipschitz smooth $f_\xi$. Since $\sum_{t=1}^T \frac{1}{t}\leq \log(T)+1$, it holds 
\begin{equation*}
\frac{1}{T}\sum_{t=1}^T\Delta_f(m_t)\leq  \frac{1}{T}\sum_{t=1}^T \frac{S_f \sigma_g^2}{2m_t}\leq \frac{S_f \sigma_g^2(\log(T)+1)}{2T}\leq \cO(\eps).
\end{equation*}
Setting $T=\tilde \cO(\eps^{-1})$, the total sample complexity is $\sum_{t=1}^T (m_t+1) = \cO(T^2)=\tilde \cO(\eps^{-2})$ for objectives with Lipschitz smooth $f_\xi$.
\end{proof}

\subsection{Global Convergence of BSGD for Convex Objectives}
We prove Theorem \ref{thm:convex}, the convex case for which   Assumption \ref{ass:convex} holds with $\mu=0$.
\begin{proof}
Recall that 
\begin{equation*}
\EE[F(\hat{x}_T)-F(x^*)] 
\leq 
\frac{1}{T}\EE \sum_{t=1}^T[F(x_t)-F(x^*)].
\end{equation*}
Since $\beta = 0$ and $\mu=0$, summing up \eqref{eq:key} from $t=1$ to $T$,
\begin{equation*}
\begin{split}
\frac{1}{T} \sum_{t=1}^T\EE[F(x_t)-F(x^*)]
\leq &
\frac{1}{T}\sum_{t=1}^T\left[ 2\Delta_f(m_t)+\frac{a_t-a_{t+1}}{\gamma_t}  + \frac{1}{2}\gamma_t\EE \norm{\nabla_{x} \hat p(x_t,\xi_t)}_2^2\right]\\
\leq & 
\frac{1}{T}\sum_{t=1}^T\left[ 2\Delta_f(m_t)+\frac{1}{2}\gamma_t M^2\right] + \frac{1}{T}\sum_{t=2}^T a_t\left(\frac{1}{\gamma_{t}} -\frac{1}{\gamma_{t-1}}\right) +\frac{1}{\gamma_{1}T}a_1.
\end{split}
\end{equation*}
Plugging constant stepsizes $\gamma_t=\gamma$ and  $a_1=\|x_1-x^*\|_2^2/2$, we have
\begin{equation*}
\EE[F(\hat{x}_T)-F(x^*)] \leq  \frac{2}{T}\sum_{t=1}^T\Delta_f(m_t) +\frac{1}{2}\gamma M^2 +\frac{\|x_1-x^*\|_2^2}{2T\gamma}.
\end{equation*}
Setting $\gamma = \frac{c}{\sqrt{T}}$, we have the desired result
\begin{equation}
\label{eq:convexcase}
\EE[F(\hat{x}_T)-F(x^*)] \leq  \frac{2}{T}\sum_{t=1}^T\Delta_f(m_t) +\frac{M^2c^2+\|x_1-x^*\|_2^2}{2c\sqrt{T}}.
\end{equation}
\end{proof}

Comparing to \citet{nemirovski2009robust} and \citet{hazan2007logarithmic}, \eqref{eq:convexcase} has an extra term $\frac{2}{T}\sum_{t=1}^T \Delta_f(m_t)$ that represents the average estimation bias of the function value estimator $\hat p(x,\xi_t)$ over $F(x)$. 

\begin{cor}
\label{cor:convex}
Under the same assumptions as Theorem \ref{thm:convex}, to achieve an $\eps$-optimal solution, 
the total sample complexity required by BSGD is $ \cO(\eps^{-4})$ for convex CSO objectives with Lipschitz continuous $f_\xi$ and $ \cO(\eps^{-3})$ for convex CSO objectives with Lipschitz smooth $f_\xi$.
\end{cor}
The sample complexity is achieved for either fixed mini-batch size $m_t = m$ or the time-varying mini-batch sizes $m_t = t$ for Lipschitz continuous $f_\xi$ or $m_t = \lceil\sqrt{t}\rceil$ for Lipschitz smooth $f_\xi$.

\begin{proof}
Let $T = \cO(\eps^{-2})$. For fixed inner batch sizes $m_t = m$, the selection of $m$ is obvious by definition of $\Delta_f(m_t)$.

For time-varying batch sizes, when $f_\xi$ is Lipschitz continuous, let $m_t = t$. Invoking $\sum_{t=1}^T 1/\sqrt{t} \leq 2\sqrt{T}$, we have,
\begin{equation*}
\EE[F(\hat{x}_T)-F(x^*)] \leq  \frac{2}{T}\sum_{t=1}^T\frac{L_f\sigma_g}{\sqrt{t}} +\frac{M^2c^2+\|x_1-x^*\|_2^2}{2c\sqrt{T}}\leq \frac{4L_f\sigma_g}{\sqrt{T}} +\frac{M^2c^2+\|x_1-x^*\|_2^2}{2c\sqrt{T}}\leq \eps.
\end{equation*}
The sample complexity is $\sum_{t=1}^T (t+1) = \cO(T^2) = \cO(\eps^{-4})$.

When $f_\xi$ is Lipschitz smooth, letting $m_t = \lceil \sqrt{t}\rceil$, we have
\begin{equation*}
\EE[F(\hat{x}_T)-F(x^*)] \leq  \frac{2}{T}\sum_{t=1}^T\frac{S_f\sigma_g^2}{2\sqrt{t}} +\frac{M^2c^2+\|x_1-x^*\|_2^2}{2c\sqrt{T}}\leq \frac{2S_f\sigma_g^2}{\sqrt{T}} +\frac{M^2c^2+\|x_1-x^*\|_2^2}{2c\sqrt{T}}\leq \eps.
\end{equation*}
The sample complexity is $\sum_{t=1}^T (\sqrt{t}+1)=\cO(T^{3/2}) = \cO(\eps^{-3})$.
\end{proof}

\subsection{Stationarity Convergence of BSGD for Weakly Convex Objectives}
	We prove Theorem \ref{thm:wc_convergence}. In this case, Assumption \ref{ass:convex} with $\mu<0$ implies that $F(x)$ is $|\mu|$-weakly convex. For simplicity, we denote $ x'_t:=\mathrm{prox}_{\lambda F}(x_t) $.   $\lambda$ is specified later in the proof.
	\begin{proof}
		By the definition of Moreau envelope, we have for any $ \hat{\mu}>|\mu| $,
		\begin{equation}
		\label{eq:nc_key}
		\begin{split}
		F_{1/\hat{\mu}}(x_{t+1})
		\leq\ &
		F(x'_t)+\frac{\hat{\mu}}{2}||x'_t-x_{t+1}||^2\\
		\leq\ &
		F(x'_t)+\hat{\mu}\gamma_t\nabla \widehat{p}(x_t,\xi_t)^\top(x'_t-x_{t+1})+\frac{\hat{\mu}}{2}||x'_t-x_t||^2-\frac{\hat{\mu}}{2}||x_{t+1}-x_t||^2\\
		=\ &
		F_{1/\hat{\mu}}(x_t)+\hat{\mu}\gamma_t\nabla \widehat{p}(x_t,\xi_t)^\top(x'_t-x_{t+1})-\frac{\hat{\mu}}{2}||x_{t+1}-x_t||^2\\
		=\ &
		F_{1/\hat{\mu}}(x_t)+\hat{\mu}\gamma_t\nabla \widehat{p}(x_t,\xi_t)^\top(x'_t-x_t)+\hat{\mu}\gamma_t\nabla \widehat{p}(x_t,\xi_t)^\top (x_t-x_{t+1})-\frac{\hat{\mu}}{2}||x_{t+1}-x_t||^2\\
		=\ &
		F_{1/\hat{\mu}}(x_t)+\hat{\mu}\gamma_t\nabla \widehat{p}(x_t,\xi_t)^\top(x'_t-x_t)+\frac{\hat{\mu}\gamma_t^2\|\nabla \widehat{p}(x_t,\xi_t)\|_2^2}{2},\\
		\end{split}
		\end{equation}
		where the second inequality comes from the triangle inequality, and the last equality comes from plugging in $x_{t+1}-x_{t}$.  By weak convexity of $ \widehat{p}(\cdot) $, we have
		\begin{equation}
		\label{eq:nc_key2}
		\begin{split}
		&\nabla\widehat{p}(x_t,\xi_t)^\top(x'_t-x_t)\\
		\leq\ &
		\widehat{p}(x'_t,\xi_t)-\widehat{p}(x_t,\xi_t)+\frac{|\mu|}{2}||x'_t-x_t||^2\\
		\leq\ &
		\underbrace{\widehat p(x'_t,\xi_t)- p(x'_t,\xi_t)}_{:=\zeta_{t1}}
		+ \underbrace{p(x'_t,\xi_t)- p(x_t,\xi_t)}_{:=\zeta_{t2}}
		+ \underbrace{p(x_t,\xi_t) - \widehat p(x_t,\xi_t)}_{:=\zeta_{t3}}
		+\frac{|\mu|}{2}||x'_t-x_t||^2.
		\end{split}
		\end{equation}
        By definition, $ \mathbb{E}_{\xi_t}p(x,\xi_t)=F(x) $. Invoking Lemma \ref{lm:bias}, $|\EE\zeta_{t1}|\leq \Delta_f(m_t)$, $|\EE\zeta_{t3}|\leq \Delta_f(m_t)$. Combining   \eqref{eq:nc_key} and \eqref{eq:nc_key2}, taking expectation over $\{\xi_t,\{\eta_{tj}\}_{j=1}^{m_t}\}$ on both sides, and using the fact that $ \EE\|\nabla\widehat{p}(x,\xi)\|_2^2 \leq M^2$, we have
        \begin{equation*}
        F_{1/\hat{\mu}}(x_{t+1})-  F_{1/\hat{\mu}}(x_t)\leq \hat{\mu}\gamma_t(2\Delta_f(m_t)+F(x'_t)-F(x_t)+\frac{|\mu|}{2}||x'_t-x_t||^2)+\frac{\hat{\mu}\gamma_t^2M^2}{2}.
        \end{equation*}
        Dividing $\hat\mu$ on both sides, rearranging and summing up from $t=1$ to $T$, we have
		\begin{equation}
		\label{eq:nonconvex_expect_1}
		\begin{split}
		& \sum_{t=1}^{T}\gamma_t\Big(F(x_t)-F(x'_t)-\frac{|\mu|}{2}||x'_t-x_t||^2\Big)\\
		\leq \ &
		\frac{1}{\hat{\mu}}\Big(F_{1/\hat{\mu}}(x_1)-F_{1/\hat{\mu}}(x_{T+1})+\frac{\hat{\mu}M^2\sum_{t=1}^{T}\gamma_t^2}{2}\Big)+2\sum_{t=1}^{T}\gamma_t\Delta_f(m_t).  
		\end{split}
		\end{equation}
		We divide $ \sum_{t=1}^{T}\gamma_t $ on both sides of the inequality above. Recall the definition of the output of the algorithm $ \hat{x}_R $. Since $\gamma_t/\sum_{t=1}^{T}\gamma_t=1/T$ due to the constant stepsize and $\hat{x}_R$ is selected from $\{x_1,...,x_T\}$ with equal probability, we have
		\begin{equation}
		\begin{split}
		 & \mathbb{E}\Big[F(\hat{x}_R)-F(\hat{x}'_R)-\frac{|\mu|}{2}||\hat{x}'_R-\hat{x}_R||^2\Big]\\
		\leq &
		\frac{F_{1/\hat{\mu}}(x_1)-F_{1/\hat{\mu}}(x_{T+1})+\frac{1}{2}\hat{\mu}M^2\sum_{t=1}^{T}\gamma_t^2+2\hat{\mu}\sum_{t=1}^{T}\gamma_t\Delta_f(m_t)}{\hat{\mu}\sum_{t=1}^{T}\gamma_t}.   
		\end{split}
		\label{eq:nonconvex_expect_2}
		\end{equation}	
		Noticing that $ F(z)+\frac{\hat{\mu}}{2}||z-x||^2 $ is $ (\hat{\mu}-|\mu|) $-strongly convex if $ \hat{\mu}>|\mu| $. Setting $\lambda = 1/\hat \mu$, we have
		\begin{equation*}
		\begin{split}
		&F(x_t)-F(x'_t)-\frac{|\mu|}{2}||x'_t-x_t||^2\\
		=\ &
		(F(x_t)+\frac{\hat{\mu}}{2}||x_t-x_t||^2)-\big(F(x'_t)+\frac{\hat{\mu}}{2}||x'_t-x_t||^2\big)+\frac{\hat{\mu}-|\mu|}{2}||x'_t-x_t||^2\\
		\geq\ &
		(\hat{\mu}-|\mu|)||x'_t-x_t||^2
		=\ 
		\frac{\hat{\mu}-|\mu|}{\hat{\mu}^2}\mathcal{G}^2_{1/\hat{\mu} F}(x_t),
		\end{split}
		\end{equation*}
		where the last inequality uses the strong convexity of $F(z)+\frac{\hat \mu}{2}\|z-x\|_2^2$.
		Recall that $\mathcal{G}_{\lambda F}(x):=\frac{1}{\lambda}||\mathrm{prox}_{\lambda F}(x)-x||_2$.
		Combining with \eqref{eq:nonconvex_expect_2}, we obtain
		\begin{equation}
		\label{eq:nonconvex_proof_core}
		\mathbb{E}\big[\mathcal{G}^2_{1/\hat{\mu} F}(\hat{x}_R)\big]
		\leq
		\frac{\hat{\mu}}{\hat{\mu}-|\mu|}\frac{F_{1/\hat{\mu}}(x_1)-F_{1/\hat{\mu}}(x_{T+1})+\frac{1}{2}\hat{\mu}M^2\sum_{t=1}^{T}\gamma_t^2+2\hat{\mu}\sum_{t=1}^{T}\gamma_t\Delta_f(m_t)}{\sum_{t=1}^{T}\gamma_t}.
		\end{equation}
		
		Plugging $\gamma_t =c/\sqrt{T}$ and $\hat \mu=2|\mu|$ into the expression above, we have
		\begin{equation*}
		\begin{split}
		& \mathbb{E}\big[\mathcal{G}^2_{1/(2|\mu|) F}(\hat{x}_R)\big]\\
		\leq\ &
		2\frac{F_{1/(2|\mu|)}(x_1)-F_{1/(2|\mu|)}(x_{T+1})+|\mu| M^2T\frac{c^2}{T}+4|\mu|\frac{c}{\sqrt{T}}\sum_{t=1}^{T}\Delta_f(m_t)}{T\cdot\frac{c}{\sqrt{T}}}\\
		=\ &
		2\frac{F_{1/(2|\mu|)}(x_1)-F_{1/(2|\mu|)}(x_{T+1})+|\mu| M^2c^2}{{c\sqrt{T}}}
		+\frac{8|\mu|\sum_{t=1}^{T}\Delta_f(m_t)}{T}.
		\end{split}
		\end{equation*} 
		By the fact that $F_{1/(2|\mu|)}(x_{T+1})\geq \inf_{x\in\mathcal{X}}F(x)$, we conclude the proof.
	\end{proof}

\begin{cor}
\label{cor:weakly_convex}
Under the same assumptions as Theorem \ref{thm:wc_convergence}, to achieve an $\eps$-stationary point, the total sample complexity required by BSGD is at most $\cO(\eps^{-8})$. If further assuming Lipschitz smooth  $f_\xi$, the sample complexity is at most $\cO(\eps^{-6})$.
\end{cor}
The proof and batch size selection are the same as the convex case. We also provide a convergence guarantee using decaying stepsizes.
\begin{cor}
\label{cor:wc_decaying}
		\textbf{(Decaying Stepsizes)}
		Let $T\geq 3$, inner batch size $m_t\equiv m$, and stepsize $\gamma_t=c/\sqrt{t}\ (t=1,\cdots,T)$ with $c>0$. If the output $ \hat{x}_R $ is chosen from $ \{x_1, \ldots, x_T\} $ with $ P(\hat{x}_R = x_i)=\gamma_i/\sum_{t=1}^{T}\gamma_t,\ (i=1,\cdots,T) $,  we have,
		\begin{equation*}
		\mathbb{E}[\mathcal{G}^2_{1/(2|\mu|)F}(\hat{x}_R)]\leq 
		\frac{2F_{1/(2|\mu|)}(x_1)-2\min_{x\in\mathcal{X}}F(x)+4|\mu| c^2 M^2\ln T}{c\sqrt{T}}+8|\mu|\Delta_f(m).
		\end{equation*}
\end{cor}
	
	\begin{proof}
	Note that the argument till \eqref{eq:nonconvex_expect_2} still applies.
		
		Plugging  $\gamma_t=\frac{c}{\sqrt{t}}$ into  \eqref{eq:nonconvex_proof_core}, we have
		\begin{equation*}
		\begin{split}
		& \mathbb{E}\big[\mathcal{G}^2_{1/(2|\mu|)F}(\hat{x}_R)\big]\\
		\leq\ &
		2\cdot\Bigg(\frac{F_{1/(2|\mu|)}(x_1)-F_{1/(2|\mu|)}(x_{T+1})+|\mu| M^2\sum_{t=1}^{T}\gamma_t^2}{\sum_{t=1}^{T}\gamma_t}+\frac{4|\mu|\sum_{t=1}^{T}\gamma_t\Delta_f(m)}{\sum_{t=1}^{T}\gamma_t}\Bigg)\\
		\leq \ & 
		\frac{2F_{1/(2|\mu|)}(x_1)-2F_{1/(2|\mu|)}(x_{T+1})+2|\mu| M^2\sum_{t=1}^{T}\gamma_t^2}{\sum_{t=1}^{T}\gamma_t}+8|\mu|\Delta_f(m).
		\end{split}
		\end{equation*} 
		Note that for $T\geq 3$
		\begin{equation*}
		\begin{split}
		\sum_{t=1}^{T}t^{-\frac{1}{2}}&\geq\int_1^{T+1}t^{-\frac{1}{2}}dt=2(\sqrt{T+1}-1)\geq\sqrt{T};\\
		\sum_{t=1}^{T}t^{-1}&\leq1+\int_2^{T}t^{-1}dt=1+\ln  T\leq2\ln T.
		\end{split}
		\end{equation*}
		We conclude the proof.
\end{proof}

\subsection{Stationary Convergence of BSpiderBoost for Nonconvex Smooth Objectives}
To analyze the convergence, define the following auxiliary functions:
\begin{equation*}
	\widehat{F}_m(x)\coloneqq\mathbb{E}_\xi\mathbb{E}_{\{\eta_i|\xi\}_{i=1}^m}f_\xi\Big(\frac{1}{m}\sum_{i=1}^{m}g_{\eta_i}(x,\xi)\Big),
\end{equation*}
where $ \{\eta_i\}_{i=1}^m $ are i.i.d samples from the conditional distribution $\PP(\eta|\xi)$.
We summarize the properties of $F$ and $\hat F_m$ as follows:
\begin{prop} 
\label{prop:BSpiderBoost}
Under Assumptions \ref{ass:general}, \ref{ass:Lipschitz}, it holds that
\begin{itemize}
    \item[(a).] $F(x)$ and $\hat F_m(x)$ are $S_F$-Lipschitz smooth where $ S_F=S_gL_f+S_fL_g^2 $.
    \item[(b).] $
	\mathbb{E}_\xi\|\nabla f_\xi(y)-\nabla \mathbb{E} f_\xi(y)\|_2^2\leq L_f^2
	$, $
	\mathbb{E}_{\xi,\eta}\|\nabla g_\eta(x,\xi)-\nabla \mathbb{E} g_\eta(x,\xi)\|_2^2\leq L_g^2, \|\nabla \hat F_m(x)\|_2^2\leq L_f^2L_g^2.
	$
	\item[(c).] By Lemma \ref{lm:bias},  $\|F(x)-\widehat{F}_m(x)\|_2^2\leq\frac{L_f^2\sigma_g^2}{m}$
	\item[(d).] By Lemma \ref{lm:bias_gradient}, $\|\nabla F(x)-\nabla \widehat{F}_m(x)\|^2_2\leq\frac{L_g^2S_f^2\sigma_g^2}{m}$.
\end{itemize}
\end{prop}
Note that there are other conditions under which these properties would hold, for instance, a natural sufficient condition to ensure the smoothness of $ F(\cdot) $ and $\hat F_m(\cdot)$ is when $ f $ or $ g $ is linear, i.e. $ S_f $ or $ S_g $ equals to zero. 
\begin{proof}
Denote $g(x,\xi)=\EE_{\eta|\xi}g_\eta(x,\xi)$, $\hat g(x,\xi) := \frac{1}{m}\sum_{j=1}^m  g_{\eta_j}(x,\xi)$.
By definition,
$$ 
\nabla F(x)
=
\nabla\mathbb{E}_{\xi}\Big[f_{\xi}\Big(g(x,\xi)\Big)\Big] 
=
\mathbb{E}_{\xi}\Big[\nabla \Big(f_{\xi}\big(g(x,\xi)\big)\Big)\Big] 
=
\mathbb{E}_{\xi}\Big[\nabla f_{\xi}\big(g(x,\xi)\big)\cdot\nabla g(x,\xi)\Big] .
$$
Note that for each fixed $ \xi $, we have
\begin{equation*}
\begin{split}
&\|\nabla \Big(f_\xi(g(x,\xi))\Big)-\nabla \Big(f_\xi(g(y,\xi))\Big)\|_2
=
\|\nabla f_\xi(g(x,\xi))\cdot\nabla g(x,\xi)-\nabla f_\xi(g(y,\xi))\cdot\nabla g(y,\xi)\|_2\\
=\ &
\|\nabla f_\xi(g(x,\xi))\cdot\nabla g(x,\xi)-\nabla f_\xi(g(x,\xi))\cdot\nabla g(y,\xi)+\nabla f_\xi(g(x,\xi))\cdot\nabla g(y,\xi)-\nabla f_\xi(g(y,\xi))\cdot\nabla g(y,\xi)\|_2\\
\leq\ &
\|\nabla f_\xi(g(x,\xi))\cdot\big(\nabla g(x,\xi)-\nabla g(y,\xi)\big)\|_2+\|\big(\nabla f_\xi(g(x,\xi))-\nabla f_\xi(g(y,\xi))\big)\cdot\nabla g(y,\xi)\|_2\\
\leq\ &
L_f\|\nabla g(x,\xi)-\nabla g(y,\xi)\|_2+S_f\|\big(g(x,\xi)-g(y,\xi)\big)\cdot\nabla g(y,\xi)\|_2\\
\leq\ &
(S_gL_f+S_fL_g^2)\|x-y\|_2,
\end{split}
\end{equation*}
where the last inequality comes from Lipschitz continuity  and Lipschitz smoothness of $g_\eta (\cdot,\xi) $. Similarly,
$$ 
\nabla \Big(f_\xi\big(\hat g(x,\xi)\big)\Big)
=
\nabla f_\xi\big(\hat g(x,\xi)\big)^\top\nabla \hat g(x,\xi).
$$
\begin{equation*}
\begin{split}
&\|\nabla \Big(f_\xi(\hat g(x,\xi))\Big)-\nabla \Big(f_\xi(\hat g(y,\xi))\Big)\|_2\\
=\ &
\|\nabla \hat g(x,\xi)^\top\nabla f_\xi(\hat g(x,\xi))-\nabla \hat g(y,\xi)^\top\nabla f_\xi(\hat g(y,\xi))\|_2\\
=\ &
\|\nabla \hat g(x,\xi)^\top\nabla f_\xi(\hat g(x,\xi))-\nabla \hat g(y,\xi)^\top\nabla f_\xi(\hat g(x,\xi))+\nabla \hat g(y,\xi)^\top\nabla f_\xi(\hat g(x,\xi))-\nabla \hat g(y,\xi)^\top\nabla f_\xi(\hat g(y,\xi))\|_2\\
\leq\ &
\|\big(\nabla \hat g(x,\xi)-\nabla \hat g(y,\xi)\big)^\top\nabla f_\xi(\hat g(x,\xi))\|_2+\|\nabla \hat g(y,\xi)^\top\big(\nabla f_\xi(\hat g(x,\xi))-\nabla f_\xi(\hat g(y,\xi))\big)\|_2\\
\leq\ &
L_f\|\nabla \hat g(x,\xi)-\nabla \hat g(y,\xi)\|_2+S_f\|\hat g(x,\xi)-\hat g(y,\xi)\|_2\cdot\|\nabla g(y,\xi)\|_2\\
\leq\ &
(S_gL_f+S_fL_g^2)\|x-y\|_2.
\end{split}
\end{equation*}
It concludes the proof of Proposition \ref{prop:BSpiderBoost}(a).

As for Proposition \ref{prop:BSpiderBoost}(b),
note that for any random variables $ X $, we have $ \mathbb{E}\|X-EX\|_2^2\leq\mathbb{E}\|X\|_2^2 $. It implies that
\begin{equation}
	\begin{split}
	\mathbb{E}_\xi\|\nabla f_\xi(y)-\nabla \mathbb{E} f_\xi(y)\|_2^2
	\leq\ &
	\mathbb{E}_\xi\|\nabla f_\xi(y)\|_2^2
	\leq
	L_f^2,\\
	\mathbb{E}_{\xi,\eta}\|\nabla g_\eta(x,\xi)-\nabla \mathbb{E} g_\eta(x,\xi)\|_2^2
	\leq\ &
	\mathbb{E}_{\xi,\eta}\|\nabla g_\eta(x,\xi)\|_2^2
	\leq
	L_g^2.
	\end{split}
\end{equation}
It further holds that
\begin{equation}
	\begin{split}
	&\mathbb{E}_\xi\mathbb{E}_{\{\eta_i|\xi\}_{i=1}^m}\Big\|\nabla\Big(f_\xi\Big(\hat g(x,\xi)\Big)\Big)-\nabla\widehat{F}_m(x)\Big\|_2^2\\
	\leq\ &
	\mathbb{E}_\xi\mathbb{E}_{\{\eta_i|\xi\}_{i=1}^m}\Big\|\nabla\Big(f_\xi\Big(\hat g(x,\xi)\Big)\Big)\Big\|_2^2\\
	=\ &
	\mathbb{E}_\xi\mathbb{E}_{\{\eta_i|\xi\}_{i=1}^m}\Big\|\nabla f_\xi\Big(\hat g(x,\xi)\Big)^\top\nabla\Big(\hat g(x,\xi)\Big)\Big\|_2^2\\
	\leq\ &
	\mathbb{E}_\xi\mathbb{E}_{\{\eta_i|\xi\}_{i=1}^m}\Big\|\nabla f_\xi\Big(\hat g(x,\xi)\Big)\Big\|_2^2\cdot\Big\|\frac{1}{m}\sum_{i=1}^{m}\nabla g_{\eta_i}(x,\xi)\Big\|_2^2\\
	\leq\ &
	L_f^2\mathbb{E}_\xi\mathbb{E}_{\{\eta_i|\xi\}_{i=1}^m}\Big[\frac{1}{m}\sum_{i=1}^{m}\Big\|\nabla g_{\eta_i}(x,\xi)\Big\|_2^2\Big]
	\leq
	L_f^2L_g^2.
	\end{split}
\end{equation}

Proposition \ref{prop:BSpiderBoost}(c) and (d) are direct implications of Lemma \ref{lm:bias} and Lemma \ref{lm:bias_gradient}, respectively.
	 
\end{proof}

Recall the gradient estimator $v_t$  of BSpiderBoost
\begin{equation}
	v_t=
	\begin{cases}
	\nabla F^{N_2}_m(x_t)-\nabla F^{N_2}_m(x_{t-1})+v_{t-1} & (n_t-1)q+1\leq t\leq n_tq-1,\\
	\nabla F^{N_1}_m(x_t) & t=(n_t-1)q,
	\end{cases}
\end{equation}
where $n_t=\lceil t/q \rceil$.
Different from a key step in the analysis of SPIDER related literature \citep{fang2018spider,wang2018spiderboost}, the sequence $ \{v_t-\nabla F(x_t)\} $ in our work is not a martingale sequence because $ \nabla F^{N_2}_m(x_t)-\nabla F^{N_2}_m(x_{t-1}) $ is not an unbiased estimator of $ \nabla F(x_t)-\nabla F(x_{t-1}) $. As a result, the error would accumulate during iterations. 
To fix the issue, we consider the following sequence in our analysis
\begin{equation}
	\{v_t-\nabla \widehat{F}_m(x_t)\}_t.
\end{equation}
which is a martingale. The key lemma of SPIDER is formulated as follows:
\begin{lm}[\cite{fang2018spider}, Lemma 1]
	\label{lm:BSpB_key_revise}
	Denote $ n_t=\lceil t/q \rceil $, then for all $ (n_t-1)q+1\leq t\leq n_tq-1 $, we have
	\begin{equation}
	\mathbb{E}\|v_t-\nabla \widehat{F}_m(x_t)\|_2^2\leq\frac{S_F^2}{N_2}\mathbb{E}\|x_t-x_{t-1}\|_2^2+\mathbb{E}\|v_{t-1}-\nabla \widehat{F}_m(x_{t-1})\|_2^2,
	\end{equation}
	and by telescoping and Proposition \ref{prop:BSpiderBoost} we have
	\begin{equation}
	\begin{split}
	\mathbb{E}\|v_t-\nabla \widehat{F}_m(x_t)\|_2^2
	\leq\ &
	\mathbb{E}\|v_{(n_t-1)q}-\nabla \widehat{F}_m(x_{(n_t-1)q})\|_2^2+\sum_{i=(n_t-1)q}^{t-1}\frac{S_F^2}{N_2}\mathbb{E}\|x_{i+1}-x_i\|_2^2\\
	\leq\ &
	\frac{L_f^2L_g^2}{N_1}+\sum_{i=(n_t-1)q}^{t-1}\frac{S_F^2}{N_2}\mathbb{E}\|x_{i+1}-x_i\|_2^2.
	\end{split}
	\end{equation}
\end{lm}

\begin{prop}\label{prop:BSpiderBoost_General_Framework_I_revise}
	For each iteration, we have
	\begin{equation}
	\begin{split}
	F(x_{t+1})
	\leq\ &
	F(x_t)-\frac{\gamma(1-S_F\gamma)}{2}\|v_t\|_2^2+\gamma\Big(\|\nabla F(x_t)-\nabla\widehat{F}_m(x_t)\|_2^2+\|\nabla\widehat{F}_m(x_t)-v_t\|_2^2\Big).
	\end{split}
	\end{equation}
\end{prop}

\begin{proof}
	By smoothness of $ F(\cdot) $, we have
	\begin{equation}
	\begin{split}
	&F(x_{t+1})
	\leq
	F(x_t)+\nabla F(x_t)^\top (x_{t+1}-x_t)+\frac{S_F}{2}||x_{t+1}-x_t||^2\\
	=\ &
	F(x_t)-\gamma\nabla F(x_t)^\top v_t+\frac{S_F\gamma^2}{2}||v_t||^2\\
	=\ &
	F(x_t)
	-\gamma
	v_t^\top(\nabla F(x_t)-\nabla \widehat{F}_m(x_t)+\nabla\widehat{F}_m(x_t)-v_t)-2\cdot\frac{\gamma}{4}\|v_t\|_2^2
	-\frac{\gamma(1-S_F\gamma)}{2}\|v_t\|_2^2,\\
	\leq\ &
	F(x_t)-\frac{\gamma(1-S_F\gamma)}{2}\|v_t\|_2^2+\gamma\Big(\|\nabla F(x_t)-\nabla\widehat{F}_m(x_t)\|_2^2+\|\nabla\widehat{F}_m(x_t)-v_t\|_2^2\Big),
	\end{split}
	\end{equation}
	the last inequality comes from Young's inequality, more precisely, i.e.
	\begin{equation}
	-\gamma
	v_t^\top(\nabla F(x_t)-\nabla \widehat{F}_m(x_t))-\frac{\gamma}{4}\|v_t\|_2^2
	\leq
	\gamma\|\nabla F(x_t)-\nabla \widehat{F}_m(x_t)\|_2^2.
	\end{equation}
\end{proof}

\begin{prop}
\label{prop:BSpiderBoost_General_Framework_II_revise}
    Denote
	$$
	\delta:=\frac{L_g^2S_f^2\sigma_g^2}{m}+\frac{L_f^2L_g^2}{N_1}
	$$
	for each iteration, we have
	\begin{equation}
	F(x_{t+1})
	\leq
	F(x_t)-\frac{\gamma(1-S_F\gamma)}{2}\|v_t\|_2^2+\gamma\delta+\gamma^3\sum_{i=(n_t-1)q}^{t-1}\frac{S_F^2}{N_2}\mathbb{E}\|v_i\|^2
	\end{equation}
\end{prop}

\begin{proof}
	Apply Proposition \ref{prop:BSpiderBoost}, \ref{prop:BSpiderBoost_General_Framework_I_revise} and Lemma \ref{lm:BSpB_key_revise}, we have
	\begin{equation}
	\begin{split}
	F(x_{t+1})
	\leq\ &
	F(x_t)-\frac{\gamma(1-S_F\gamma)}{2}\|v_t\|_2^2+\gamma\Big(\|\nabla F(x_t)-\nabla\widehat{F}_m(x_t)\|_2^2+\|\nabla\widehat{F}_m(x_t)-v_t\|_2^2\Big)\\
	\leq\ &
	F(x_t)-\frac{\gamma(1-S_F\gamma)}{2}\|v_t\|_2^2+\gamma\Big(\frac{L_f^2\sigma_{gg}^2+L_g^2S_f^2\sigma_g^2}{m}+\frac{L_f^2L_g^2}{N_1}+\sum_{i=(n_t-1)q}^{t-1}\frac{S_F^2}{N_2}\mathbb{E}\|x_{i+1}-x_i\|_2^2\Big)\\
	=\ &
	F(x_t)-\frac{\gamma(1-S_F\gamma)}{2}\|v_t\|_2^2+\gamma^3\sum_{i=(n_t-1)q}^{t-1}\frac{S_F^2}{N_2}\mathbb{E}\|v_i\|^2+\gamma\delta.
	\end{split}
	\end{equation}
	which concludes the proof.
\end{proof}

Now we proceed to the proof of the convergence rate of BSpiderBoost.
\begin{proof}[Proof of Theorem \ref{thm:BSpiderBoost_CSO_revise}]
	Based on Proposition \ref{prop:BSpiderBoost_General_Framework_II_revise}, we telescope in each epoch and take expectation, 
	\begin{equation}
	\begin{split}
	&\mathbb{E}\big[F(x_{t+1})\big]\\
	\leq\ &
	\mathbb{E}\Big[F(x_{(n_t-1)q})-\frac{\gamma(1-S_F\gamma)}{2}\sum_{j=(n_t-1)q}^{t}\|v_j\|_2^2+\gamma^3\sum_{j=(n_t-1)q}^{t}\sum_{i=(n_t-1)q}^{j-1}\frac{S_F^2}{N_2}\mathbb{E}\|v_i\|^2+\gamma\sum_{j=(n_t-1)q}^{t}\delta\Big]\\
	\leq\ &
	\mathbb{E}\Big[F(x_{(n_t-1)q})-\frac{\gamma(1-S_F\gamma)}{2}\sum_{j=(n_t-1)q}^{t}\|v_j\|_2^2+\gamma^3\sum_{j=(n_t-1)q}^{n_tq-1}\sum_{i=(n_t-1)q}^{t}\frac{S_F^2}{N_2}\mathbb{E}\|v_i\|^2+\gamma\sum_{j=(n_t-1)q}^{t}\delta\Big]\\
	=\ &
	\mathbb{E}\Big[F(x_{(n_t-1)q})-\frac{\gamma(1-S_F\gamma)}{2}\sum_{j=(n_t-1)q}^{t}\|v_j\|_2^2+\frac{\gamma^3S_F^2q}{N_2}\sum_{i=(n_t-1)q}^{t}\mathbb{E}\|v_i\|^2+\gamma\sum_{j=(n_t-1)q}^{t}\delta\Big]\\	
	\leq\ &
	\mathbb{E}\Big[F(x_{(n_t-1)q})-\sum_{i=(n_t-1)q}^{t}\Big(\beta\|v_i\|_2^2-\gamma\delta\Big)\Big],
	\end{split}
	\end{equation}
	the second inequality follows from $ n_tq>t $ and $ j-1<t $. Further telescoping for all iterations, we have
	\begin{equation}
	\begin{split}
	&\mathbb{E}\big[F(x_T)-F(x_0)\big]
	\leq
	\mathbb{E}\Big[-\sum_{i=0}^{T-1}\Big(\beta\|v_i\|_2^2-\gamma\delta\Big)\Big]
	=
	\gamma T\delta-\beta\sum_{i=0}^{T-1}\mathbb{E}\|v_i\|_2^2.
	\end{split}
	\end{equation}
	As a result, for $ x_S $, its corresponding $ v_S $ satisfies that
	\begin{equation}\label{eq:BSpiderBoost_vs_revise}
	\mathbb{E}\|v_S\|_2^2
	\leq
	\frac{1}{T}\sum_{i=0}^{T-1}\mathbb{E}\|v_i\|_2^2
	\leq
	\frac{F(x_0)-\mathbb{E}F(x_T)}{\beta T}+\frac{\gamma}{\beta}\delta
	\leq
	\frac{\Delta}{\beta T}+\frac{\gamma}{\beta}\delta.
	\end{equation}
    By substituting the parameter settings, we could further show that
	\begin{equation}
	\begin{split}
	\mathbb{E}\Big[\frac{1}{3}\|\nabla F(x_S)\|_2^2\Big]
	\leq\ &
	\mathbb{E}\Big[\|\nabla F(x_S)-\nabla\widehat{F}_m(x_S)\|_2^2+\|\nabla\widehat{F}_m(x_S)-v_S\|_2^2+\|v_S\|_2^2\Big]\\
	\leq\ &
	\frac{L_g^2S_f^2\sigma_g^2}{m}
	+
	\Big(\frac{L_f^2L_g^2}{N_1}+\mathbb{E}\sum_{i=(n_S-1)q}^{S-1}\frac{S_F^2}{N_2}\mathbb{E}\|x_{i+1}-x_i\|_2^2\Big)
	+
	\Big(\frac{\Delta}{\beta T}+\frac{\gamma}{\beta}\delta\Big)\\
	\leq\ &
	\delta+\frac{\Delta}{\beta T}+\frac{\gamma}{\beta}\delta+\frac{S_F^2}{N_2}\mathbb{E}\sum_{i=(n_S-1)q}^{S-1}\mathbb{E}\|x_{i+1}-x_i\|_2^2\\
	\leq\ &
	\delta+\frac{\Delta}{\beta T}+\frac{\gamma}{\beta}\delta+\frac{\gamma^2S_F^2}{N_2}\cdot\frac{q}{T}\sum_{i=0}^{T-1}\mathbb{E}\|v_i\|_2^2\\
	\leq\ &
	\delta+\frac{\Delta}{\beta T}+\frac{\gamma}{\beta}\delta+\frac{\gamma^2S_F^2q}{N_2}\Big(\frac{\Delta}{\beta T}+\frac{\gamma}{\beta}\delta\Big)\\
	=\ &
	\Big(1+\frac{\gamma}{\beta}+\frac{\gamma^3S_F^2q}{N_2\beta}\Big)\delta+\frac{1}{\beta T}\Big(1+\frac{\gamma^2S_F^2q}{N_2}\Big)\Delta\\
	\leq\ &
	\Big(1+\frac{1}{2\beta S_F}+\frac{1}{16\beta S_F}\Big)\delta
	+
	\frac{1}{\beta T}\Big(1+\frac{1}{8}\Big)\Delta
	\end{split}
	\end{equation}
	where the second inequality comes from Proposition \ref{prop:BSpiderBoost}(d), Lemma \ref{lm:BSpB_key_revise}, and Equation \eqref{eq:BSpiderBoost_vs_revise}; the third inequality follows from the fact that the probability such that $ n_S=1 $ or $ 2,\cdots,n_T $ is less or equal to $ q/T $; recall the definition of $\delta$ and multiply both sides by 3, then we have
	\begin{equation}
	    \begin{split}
	        \mathbb{E}\|\nabla F(x_S)\|_2^2
    	    \leq\ &
    	    \Big(3+\frac{3}{2\beta S_F}+\frac{3}{16\beta S_F}\Big)\Big(\frac{L_g^2S_f^2\sigma_g^2}{m}+\frac{L_f^2L_g^2}{N_1}\Big)+\frac{4\Delta}{\beta T}\\
        	\leq\ &
        	\frac{\epsilon^2}{4}+\frac{\epsilon^2}{4}+\frac{\epsilon^2}{2}=\epsilon^2,
	    \end{split}
	\end{equation}
	By Jensen's inequality for the function $x^2$, it holds that
	\begin{equation}
	\big(\mathbb{E}\|\nabla F(x_S)\|_2\big)^2\leq \mathbb{E}\|\nabla F(x_S)\|_2^2\leq \epsilon^2.
	\end{equation}
	So $x_S$ is the stationary point we desire. The corresponding iteration complexity is 
	\begin{equation}
	\lceil T/q\rceil N_1+2TN_2=\mathcal{O}(\epsilon^{-3}),
	\end{equation}	
	and the sample complexity is 
	\begin{equation}
	\lceil T/q\rceil N_1m+2TN_2m=\mathcal{O}(\epsilon^{-5}).
	\end{equation}
	which concludes the proof.
\end{proof}

\section{Lower Bounds, Proof of Theorem \ref{thm:lower_bound_CSO}}
\label{secapp:lower}
\begin{proof}
We prove the lower bounds of $\Delta(\mathcal{A},\mathcal{F},\Phi)$ by constructing a hard instance of function and a hard instance of an oracle. 

The hard instance $F(x)$ is such that $F(x)$ is a CSO objective and it satisfies that $x=(y,z)$ and $F(x) = F^\mathrm{c}(y)+F^\mathrm{sm}(z)$. It means that $F$ is separable on $y$ and $z$. 

The oracle $\phi\in\Phi_m$ returns a biased function value and a gradient estimator of $F$. Since $F$ is separable, we specifically consider a hard oracle instance $\phi$ such that $\phi$ returns a biased function value estimator and a biased gradient estimator of $F^\mathrm{sm}(z)$ due to the compositional structure on $z$ and it returns an unbiased function value estimator and an unbiased gradient estimator of $F^\mathrm{c}(y)$. 

Based on this specific hard instance construction on the function and the oracle, we could decompose the lower bounds into two parts: one on coordinate $y$ and the other on coordinate $z$. Note that the part on coordinate $y$ is the classical lower bounds using the stochastic first-order oracle. The part on coordinate $z$ is related to the extra bias term introduced by the biased oracle.

We first consider the (strongly) convex CSO function class.
By Yao's principle~\citep{yao1977}, we have
\begin{equation}
\label{eq:lower_decomposition}
\begin{split}
& \Delta_T^*(\mathcal{A},\mathcal{F},\Phi)  \\
= &
\inf_{A \in \mathcal{A}} \sup_{F \in \mathcal{F}} \sup_{\phi\in\Phi} \EE \Delta_T^A(F,\phi,\mathcal{X}) \\
\geq &
\inf_{A \in \mathcal{A}^d} \EE_{\{V,\phi\}} \Delta_T^A(F_{V},\phi,\mathcal{X})\\
= & \inf_{A \in \mathcal{A}^d} \EE_{\{V ,\phi \}}\left[F_{V }(x_T^A(\phi ))-F_{V }(x_V^*)\right],\\
\geq & \inf_{A \in \mathcal{A}^d} \EE_{\{V,\phi\}}\left[F^\mathrm{c}_{V }(y_T^A(\phi ))-F^\mathrm{c}_{V }(y^*_V)+F^\mathrm{sm}_{V }(z_T^A(\phi ))-F^\mathrm{sm}_{V }(z_V^*)\right]   \\
\geq & 
\inf_{A \in \mathcal{A}^d} \EE_{\{V,\phi\}}\left[F^\mathrm{c}_{V }(y_T^A(\phi ))-F^\mathrm{c}_{V }(y^*_V)\right] +  \inf_{A \in \mathcal{A}^d}  \EE_{\{V,\phi\}}\left[F^\mathrm{sm}_{V }(z_T^A(\phi ))-F^\mathrm{sm}_{V }(z_V^*)\right]
\end{split}
\end{equation}
where $\cA^d$ represents the class of all deterministic algorithms, $x_T^A(\phi) = (y_T^A(\phi),z_T^A(\phi))$, $x_V^* = (y_V^*,z_V^*)$ is the minimizer of $F_V(x)$. We first consider the lower bounds on the $y$ part where the oracle returns an unbiased function value and gradient estimator. To lower bound $\inf_{A \in \mathcal{A}^d} \EE_{\{V,\phi\}}\left[F^\mathrm{c}_{V }(y_T^A(\phi ))-F^\mathrm{c}_{V }(y^*_V)\right]$, the hard instance construction of $F^c_V(y)$ will satisfy the following conditions.

\textbf{Condition I}:
\begin{itemize}[leftmargin=2em]
    \item For $V=+1$ or $V= -1$, if $yV\leq 0$, then $F^\mathrm{c}_{V}(y)\geq F^\mathrm{c}_{V}(0)$.
    \item For $V=+1$ or $V= -1$, $F^\mathrm{c}_{V}(0)-\inf_{y}F^\mathrm{c}_{V}(y)=c_0$ where $c_0\geq0$ is a constant.
\end{itemize}
Suppose that Condition I holds, we have
\begin{equation}
\label{eq:classical_part}
\begin{split}
&\inf_{A \in \mathcal{A}^d} \EE_{\{V,\phi\}}\left[F^\mathrm{c}_{V }(y_T^A(\phi ))-F^\mathrm{c}_{V }(y^*_V)\right]\\
\geq & \inf_{A \in \mathcal{A}^d} \EE_{\{V,\phi\}}\left[ (F^\mathrm{c}_{V}(0)-F^\mathrm{c}_{V }(y^*_V))\mathbb{I}\{y_T^A(\phi)V\leq 0\}\right]\\
=& \inf_{A \in \mathcal{A}^d} c_0\EE_{\{V,\phi\}}\left[ \mathbb{I}\{y_T^A(\phi)V\leq 0\}\right]\\
= & \inf_{A \in \mathcal{A}^d} c_0\PP(y_T^A(\phi)V\leq 0).
\end{split}
\end{equation}
Notice that  \eqref{eq:classical_part} requires to lower bound a probability
$$
\PP_V\{y_T^A(\phi )V\leq 0\}.
$$
For a constant $v=+1$ or $v= -1$, let $\PP^v$ denote the probability distribution of the following trajectory
$$
(y_0^A, G_v(y_0^A,\zeta_0), y_1^A(\phi),..., G_v(y_{T-1}^A(\phi),\zeta_{T-1}), y_T^A(\phi) ).
$$
It holds that 
$$
\PP_V\{y_T^A(\phi )V\leq 0\}\geq 1- \|\PP^{+1}-\PP^{-1}\|_\mathrm{TV} \geq 1-\sqrt{0.5D_\mathrm{KL}(\PP^{+1}\mid\mid \PP^{-1})},
$$
where $\|\cdot\|_\mathrm{TV}$ denotes the total variation distance of two probability distributions, $D_\mathrm{KL}$ denotes the KL divergence of two probability distributions. The first inequality holds by definition and the second inequality comes from Pinsker's inequality~\citep{cover1999elements}. Since $A\in\mathcal{A}^d$ is a deterministic algorithm, conditioned on oracle return $\phi(y_j,\zeta_j)$, $y_j^A(\phi)$ is deterministic. Conditioned on $y_j^A(\phi)$, the randomness in $\phi(y_{j}^A(\phi),\zeta_j)$ only comes from $\zeta_j$. By the chain rule of KL divergence, we have
$$
D_\mathrm{KL}\Big(\PP^{+1}\mid\mid \PP^{-1}\Big)= \sum_{t=0}^{T-1} D_\mathrm{KL}(G_{+1}(y_t^A(\phi_{+1}),\zeta_t)|y_t^A(\phi_{+1})\mid\mid G_{-1}(y_t^A\phi_{-1},\zeta_t)|y_t^A(\phi_{-1})).
$$
In our hard oracle construction, we require the following conditions:

\noindent\textbf{Condition II}
\begin{itemize}[leftmargin=2em]
    \item The gradient estimator returned by the oracle, conditioned on the query point, is a normal random variable such that conditioned on $w_t$, it holds:
$$
G_{V}(y_t^A(\phi),\zeta_t))|y_t^A(\phi)\sim\NN(\mu_V^A,\sigma^2)
$$
where $\mu_V^A$ depends on the algorithm $A$, $\sigma^2$ is the variance parameter of the oracle.
\item There exists a constant $b_y$ such that $|\mu_{+1}^A-\mu_{-1}^A|\leq b_y$.
\end{itemize}

Since the KL divergence between two normal random variables with the same variance $\sigma^2$ is known to be $\frac{(\mu_1-\mu_2)^2}{2\sigma^2}$, where $\mu_1$ and $\mu_2$ are the expectations of the two normal random variables, respectively. As a result, we have
$$
D_\mathrm{KL}(\PP^{+1}\mid\mid \PP^{-1})  = \frac{T(\mu_{+1}^A-\mu_{-1}^A)^2}{2\sigma^2}\leq \frac{Tb_y^2}{2\sigma^2}.
$$
Thus, it implies that
\begin{equation}
\label{eq:KL_divergence}
\PP_V\{y_T^A(\phi )V\leq 0\}\geq 1-\sqrt{\frac{T^2b_y^2}{4\sigma^2}}.
\end{equation}
Therefore we have the lower bounds on the first term on the right hand side of \eqref{eq:lower_decomposition}:
\begin{equation}
\label{eq:classical_final}
\inf_{A \in \mathcal{A}^d} \EE_{\{V,\phi\}}\left[F^\mathrm{c}_{V }(y_T^A(\phi ))-F^\mathrm{c}_{V }(y^*_V)\right]\geq  c_0\bigg(1-\sqrt{\frac{T^2b_y^2}{4\sigma^2}}\bigg).
\end{equation}

Now we lower bound the second term $\inf_{A \in \mathcal{A}^d}  \EE_{\{V,\phi\}}\left[F^\mathrm{sm}_{V }(z_T^A(\phi ))-F^\mathrm{sm}_{V }(z_V^*)\right]$ on the right hand side of \eqref{eq:lower_decomposition}. 
Let $\hat F^\mathrm{sm}_{m,V }(z)$ denote the expectation of the  function value estimator returned by the oracle on $z$ part. Suppose the hard instance construction satisfies the following conditions.

\textbf{Condition III}:
\begin{itemize}[leftmargin=2em]
    \item For $V=+1$ or $V= -1$, if $zV\leq 0$, then $F^\mathrm{sm}_{V}(z)\geq F^\mathrm{sm}_{V}(0)$.
    \item For $V=+1$ or $V= -1$, $F^\mathrm{sm}_{V}(0)-\hat F_{m,V}^\mathrm{sm}(z_V^*) = c_m$ where $c_m\geq0$ for any $m\geq m_0$ .
    \item For $V=+1$ or $V=-1$, $\hat F^\mathrm{sm}_{m,V}(z^*_V)- F^\mathrm{sm}_{V}(z^*_V)=c_m^\prime$, where $c_m^\prime\geq 0$ is a constant.
\end{itemize}
Condition III guarantees that for any $z$ such that $zV\leq 0$, it holds
$$
F_V^\mathrm{sm}(z) - \hat F_{m,V}^\mathrm{sm}(z_V^*) \geq F_V^\mathrm{sm}(0) - \hat F_{m,V}^\mathrm{sm}(z_V^*) \geq c_m\geq 0.
$$
We further have that
\begin{equation}
\label{eq:bias_part}
\begin{split}
&\inf_{A \in \mathcal{A}^d}  \EE_{\{V,\phi\}}\left[F^\mathrm{sm}_{V }(z_T^A(\phi ))-F^\mathrm{sm}_{V }(z_V^*)\right]\\
\geq & \inf_{A \in \mathcal{A}^d} \EE_{\{V,\phi\}} (F^\mathrm{sm}_V(0) -  F^\mathrm{sm}_{V}(z^*_V))\mathbb{I}\{z_T^A(\phi ) V\leq 0\} \\
= & \inf_{A \in \mathcal{A}^d} \EE_{\{V,\phi\}} (F^\mathrm{sm}_V(0) -  \hat F^\mathrm{sm}_{m,V}(z^*_V)+\hat F^\mathrm{sm}_{m,V}(z^*_V)- F^\mathrm{sm}_{V}(z^*_V))\mathbb{I}\{z_T^A(\phi ) V\leq 0\}\\
= &\inf_{A \in \mathcal{A}^d} (c_m+c_m^\prime)\PP\{z_T^A(\phi ) V\leq 0\}\\
\geq & \inf_{A \in \mathcal{A}^d} c_m^\prime\PP\{z_T^A(\phi ) V\leq 0\}.
\end{split}
\end{equation}
Similar as \eqref{eq:KL_divergence}, we have 
$$
\PP\{z_T^A(\phi ) V\leq 0\}\geq 1-\sqrt{\frac{T^2b_z^2}{4\sigma^2}},
$$
where $b_z\geq |\mu_1^A-\mu_{-1}^A|$ with $\mu_V^A$ as the expectation of the gradient estimator $G_V$ on $z$.
To summarize, if Condition I, II, and III hold, we have
\begin{equation}
\Delta_T^*(\mathcal{A},\mathcal{F},\Phi)\geq c_0\bigg(1-\sqrt{\frac{T^2b_y^2}{4\sigma^2}}\bigg)+ c_m^\prime\Big(1-\sqrt{\frac{T^2b_z^2}{4\sigma^2}}\Big).
\end{equation}
It remains to construct specific hard instances satisfying these conditions and demonstrate the lower bounds.
 
\paragraph{Strongly convex CSO objective with smooth outer function}
The hard instance is
\begin{equation}
\label{eq:hard_example_sc_smooth}
F_V(x) =\EE_\xi[f_\xi(\EE_{\eta|\xi} g_\eta(x,\xi))]=\frac{1}{2} y^2- V\alpha y+\frac{1}{2} z^2- V\beta z.
\end{equation}
Note that the outer function $f_\xi(y,z) = \frac{1}{2} y^2- V\alpha y+\frac{1}{2} z^2- V\beta z$ and the inner function $g_\eta(x,\xi) = (y, \eta z)$, where $\eta\sim \NN(1,1/\beta^2)$ for any $\xi$.
The minimizer of $F_V(x)$ is $  x_V^* = (V\alpha,V\beta)$. Further $F(x)$ can be decomposed as 
$ F^\mathrm{c}_V(y)=\frac{1}{2} y^2- V\alpha y$ and $F^\mathrm{sm}_V(z)=\EE_\xi[\frac{1}{2}(\EE_{\eta|\xi} \eta z)^2- V (\EE_{\eta|\xi} \eta z\beta)]=\frac{1}{2}z^2-V\beta z$.
The oracle has access to an approximation function
$$
\hat F_{m,V}(x)=F^\mathrm{c}_V(y)+\hat F^\mathrm{sm}_{m,V}(z)= \frac{1}{2} y^2- V\alpha y+\frac{1}{2} \EE \hat \eta^2 z^2- V\beta\EE \hat \eta z =\frac{1}{2} y^2- V\alpha y+\frac{1}{2} \big(1+\frac{1}{m\beta^2}\big)z^2- V\beta z,
$$
where $\hat \eta =\frac{1}{m}\sum_{j=1}^m \eta_j$, and $$
F^\mathrm{sm}_V(z) = \frac{1}{2}z^2-Vz\beta, \quad \hat F^\mathrm{sm}_V(z) = \frac{1}{2}\Big(1+\frac{\sigma_\eta^2}{m}\Big)z^2-Vz\beta.
$$.
The gradient estimator returned by the oracle  is
$$
G_V(x) =  (\nabla F^\mathrm{c}(y), \nabla \hat F_{m,V}^\mathrm{sm}(z))+\zeta),
$$
where $\zeta\sim\NN(0,\sigma^2)$. Thus $\EE G_V(x) = (\nabla F^\mathrm{c}(y), \nabla \hat F_{m,V}^\mathrm{sm}(z))) = \nabla \hat F_{m,V}(x)$. We  verify that 
\begin{itemize}[leftmargin=2em]
\item $F_V(x)$ is strongly convex and  the outer function $f_\xi(y,z)$ is $\frac{1}{2}$-Lipschitz smooth.
\item Condition I and III  hold:
\begin{itemize}
    \item $c_0 = \frac{\alpha^2}{2}$. 
    \item $c_m = \frac{\beta^2}{2}(1-\frac{1}{m\beta^2})$. Thus $c_m\geq 0$ if $m\geq 1/\beta^2$.
    \item $c_m^\prime = \frac{1}{2m}$.
\end{itemize}
\item Condition II holds:
\begin{itemize}
 \item $y-V\alpha+\zeta$ and $z-V\beta+\zeta$ are normal random variables conditioned on $V$ and $x$.
    \item $b_y = 2\alpha$, $b_z = 2\beta$.
\end{itemize}
\end{itemize}
As a result we have
\begin{equation*}
\Delta_T^*(\mathcal{A}, \mathcal{F}^+_\mathrm{CSO},\Phi_\mathrm{m} ) 
\geq 
\frac{\alpha^2}{2}(1-\sqrt{\frac{T\alpha^2}{\sigma^2}}) + \frac{1}{2m}(1-\sqrt{\frac{T\beta^2}{\sigma^2}}).
\end{equation*}
Thus there exists a hard instance with $\alpha = \frac{2}{3}\sqrt{\frac{\sigma^2}{T}}$, $\beta = \frac{1}{2}\sqrt{\frac{\sigma^2}{T}}$ such that
\begin{equation}
\label{eq:lb_sc_sm_result}
\Delta_T^*(\mathcal{A}, \mathcal{F}^+_\mathrm{CSO},\Phi_\mathrm{m} )   \geq \frac{4\sigma^2}{27T}+\frac{1}{4m},
\end{equation}
for strongly convex CSO objective with smooth outer function. As a result, to achieve $\eps$-optimality, the number of iterations should be at least $T=\cO(\eps^{-1})$, the inner batch size should be at least $m =\cO(\eps^{-1})$.

\paragraph{Strongly convex CSO objective with non-smooth outer function}
We consider the hard instance construction:
$$
F_V(x) =\EE_\xi f_\xi(\EE_{\eta|\xi} g_\eta(x,\xi))= \frac{1}{2}y^2-V\alpha y + \beta|z-V|+\beta^2(\frac{1}{2}z^2-V z),
$$
with the outer function $f_\xi(y,z) = \frac{1}{2} y^2- V\alpha y+\beta|z-V|+\beta^2(\frac{1}{2}z^2-Vz)$ and the inner function  $g_\eta(x,\xi) = (y, \eta z)$, where $\eta\sim \NN(1,1/\beta^2)$ for any $\xi$. The minimizer of $F_V(x)$ is $ x_V^* = (V\alpha, V)$. Correspondingly, we have $ F^\mathrm{c}_V(y)=\frac{1}{2} y^2- V\alpha y$ and $F^\mathrm{sm}_V(z)=\beta|z-V|+\beta ^2(\frac{1}{2}z^2-Vz)$.
The oracle has access to an approximation function
$$
\hat F_{m,V}(x) = F^\mathrm{c}_V(y)+\hat F^\mathrm{sm}_{m,V}(z)=  \frac{1}{2}y^2-V\alpha y+ \beta\EE [|\hat \eta z- V|+\beta^2(\frac{1}{2}(\hat \eta z)^2-V (\hat \eta z))].
$$
The subgradient estimator returned by the oracle  for $\hat F_{m,V}(x)$ is
$$
G_V(x) = (\nabla F^\mathrm{c}(y)+\zeta,  \nabla\hat F^\mathrm{sm}_{m,V}(z)+\zeta)
$$
where $\zeta\sim\NN(0,\sigma^2)$. We  verify that
\begin{itemize}[leftmargin=2em]
    \item $F_V(x)$ is strongly convex and the outer function $f_\xi(y,z)$ is non-smooth. 
    \item Condition I and III hold 
    \begin{itemize}
            \item $c_0 = \frac{\alpha^2}{2}$.
    \item $c_m = \beta (1-\EE |\hat \eta-1|)+\frac{\beta^2}{2}(1-\EE (\hat \eta-1)^2)$. Since $\EE |\hat \eta-1|=\frac{1}{\beta\sqrt{m}}\sqrt{\frac{2}{\pi}}$, and $\EE (\hat \eta-1)^2=\frac{1}{\beta^2 m}$, $c_m\geq 0$ if $1\geq \frac{1}{\beta\sqrt{m}}\sqrt{\frac{2}{\pi}}$ and $1\geq \frac{1}{\beta^2 m}$.
    \item $c_m^\prime = \frac{1}{\sqrt{m}}\sqrt{\frac{2}{\pi}}+\frac{1}{2m}$.
    \end{itemize}
    \item Condition II is satisfied:
    \begin{itemize}
    \item $\nabla F^\mathrm{c}(y)+\zeta$ and $\nabla\hat F^\mathrm{sm}_{m,V}(z)+\zeta$ are normal random variables conditioned on $V$ and $x$.
    \item $b_y = 2\alpha$, $b_z = 2\beta+2\beta^2$ as $|\nabla\hat F^\mathrm{sm}_{m,-}(z)-\nabla\hat F^\mathrm{sm}_{m,+}(z)|\leq 2\beta+2\beta^2$.
    \end{itemize}
\end{itemize}

As a result, we have
\begin{equation*}
\begin{split}
\Delta_T^*(\mathcal{F}^+_\mathrm{CSO}, \Phi_\mathrm{m}, \mathcal{A}) 
\geq & \frac{\alpha^2}{2}\Big(1-\sqrt{\frac{T\alpha^2}{\sigma^2}}\Big) +(\frac{1}{\sqrt{m}}\sqrt{\frac{2}{\pi}}+\frac{1}{2m})\Big(1-\sqrt{\frac{T(\beta+\beta^2)^2}{\sigma^2}}\Big)
\end{split}
\end{equation*}
As a result, there exists a hard instance with $\alpha = \frac{2}{3}\sqrt{\frac{\sigma^2}{T}}$, $\beta+\beta^2=\frac{1}{2}\sqrt{\frac{\sigma^2}{T}}$, such that
\begin{equation}
\label{eq:lb_sc_nonsmooth_result}
\Delta_T^*( \mathcal{A}, \mathcal{F}^+_\mathrm{CSO}, \Phi_\mathrm{m}) 
\geq  
\frac{4\sigma^2}{27T} + \frac{1}{2\sqrt{m}}\sqrt{\frac{2}{\pi}}+\frac{1}{4m},
\end{equation}
for strongly convex CSO objective with non-smooth outer function.

\paragraph{Convex CSO objective with smooth outer function} 
We consider a hard instance that differs from the hard instance in the strongly convex objectives with a smooth outer function only in $F_V^\mathrm{c}(y)$. 
$$
F_V(x) =\EE_\xi[f_\xi(\EE_{\eta|\xi} g_\eta(x,\xi))]=\alpha \mathbb{I}\{|y-V|>r\} (|y-V|-\frac{1}{2}r)+\alpha\mathbb{I}\{|y-V|\leq r\}\frac{1}{2r}(y-V)^2+\frac{1}{2} z^2- V z\beta.
$$
where the inner function  $g_\eta(x,\xi) = (y, \eta z)$, where $\eta\sim \NN(1,1/\beta^2)$ for any $\xi$ and the outer function 
$$f_\xi(y,z) = \alpha \mathbb{I}\{|y-V|>r\} (|y-V|-\frac{1}{2}r)+\alpha\mathbb{I}\{|y-V|\leq r\}\frac{1}{2r}(y-V)^2+\frac{1}{2} z^2- Vz\beta.$$ 
Thus $F_V(x)$ is a convex CSO objective with smooth outer function. Specifically, we have
$$
F^\mathrm{c}_V(y)=\alpha \mathbb{I}\{|y-V|>r\} (|y-V|-\frac{1}{2}r)+\alpha\mathbb{I}\{|y-V|\leq r\}\frac{1}{2r}(y-V)^2,
$$
where $r>0$. 
We construct an oracle with gradient estimator such that:
$$
G_V(x) = (\nabla F^\mathrm{c}_V(y)+ \zeta, \nabla \hat F^\mathrm{sm}_V(z) +\zeta),
$$
where $\zeta\sim\NN(0,\zeta)$.
We verify that
\begin{itemize}[leftmargin=2em]
  \item $F_V(x)$ is convex and the outer function $f_\xi(y,z)$ is $\frac{1}{r}$-Lipschitz smooth. 
    \item Condition I and III hold:
    \begin{itemize}
        \item when $0<r<1$ and $m\geq 1/\beta^2$, we have $c_0 = \alpha(1-r/2)\geq 0$, $c_m = \frac{\beta^2}{2}(1-\frac{1}{\beta^2 m})\geq 0$.
        \item $c_m^\prime = \frac{1}{2m}$.
    \end{itemize}
    \item Condition II holds:
    \begin{itemize}
     \item $\nabla F^\mathrm{c}_V(y)+ \zeta$ and $\nabla \hat F^\mathrm{sm}_V(z) +\zeta$ are a normal random variables conditioned on $V$ and $x$.
        \item $b_y = 2\alpha$, $b_z=2\beta$.
    \end{itemize}
   
\end{itemize}
As a result, we have,
\begin{equation}
\Delta_T^*(\mathcal{F}^0_\mathrm{CSO}, \Phi_\mathrm{m}, \mathcal{A})   \geq \frac{\alpha}{2}\Big(1-\sqrt{\frac{T\alpha^2}{\sigma^2}}\Big)+\frac{1}{2m}\Big(1-\sqrt{\frac{T\beta^2}{\sigma^2}}\Big).
\end{equation}
Thus there exist a hard instance with $\alpha = \frac{1}{2}\sqrt{\frac{\sigma^2}{T}}$ and $\beta = \frac{1}{2}\sqrt{\frac{\sigma^2}{T}}$ such that
\begin{equation}
\Delta_T^*(\mathcal{F}^0_\mathrm{CSO}, \Phi_\mathrm{m}, \mathcal{A})   \geq \sqrt{\frac{\sigma^2}{64T}}+ \frac{1}{4m}
\end{equation}
for convex CSO objective with smooth outer function.

\paragraph{Convex CSO objective with non-smooth outer function}
We consider the hard instance construction such that 
$$
F_V(x) =\EE_\xi f_\xi(\EE_{\eta|\xi} g_\eta(x,\xi))= \alpha|y-V| + \beta|z-V|.
$$
where the inner function is $g_\eta(x,\xi) = (y, z-\eta )$ with $\eta\sim \NN(0,1/\beta^2)$ for any $\xi$ and the outer function $f_\xi(y,z) = \alpha|y-V|+\beta|z-V|$. 
The minimizer of $F_V$ is $ x_V^* = (V, V)$.
The oracle has access to an approximation function:
$$
\hat F_{m,V}(x) = F^\mathrm{c}_V(y)+\hat F^\mathrm{sm}_{m,V}(z)=  \alpha|y-V|+\beta \EE |z-V-\hat \eta|.
$$
The subgradient estimator return by the oracle  for $\hat F_{m,V}(x)$ is
$$
G_V(x) = (\nabla F^\mathrm{c}(y)+\zeta,  \nabla\hat F^\mathrm{sm}_{m,V}(z)+\zeta)
$$
where $\zeta\sim\NN(0,\sigma^2)$ and we abuse the use of $\nabla$ to denote subgradient. We  verify that
\begin{itemize}[leftmargin=2em]
    \item $F_V(x)$ is convex and the outer function $f_\xi(y,z)$ is non-smooth. 
    \item Condition I and III are satisfied by our construction:
    \begin{itemize}
     \item $c_0 = \alpha$.
    \item $c_m = \beta(1-\EE |\hat \eta|)$, since $\EE |\hat \eta|=\frac{1}{\beta\sqrt{m}}\sqrt{\frac{2}{\pi}}$, $c_m\geq 0$ if $1\geq \frac{\sigma_\eta}{\sqrt{m}}\sqrt{\frac{2}{\pi}}$.
    \item $c_m^\prime =\beta \EE |\hat \eta| = \frac{1}{\sqrt{m}}\sqrt{\frac{2}{\pi}}$.
    \end{itemize}
    \item Condition II is satisfied:
    \begin{itemize}
    \item $\nabla F^\mathrm{c}_V(y)+\zeta$ and$ \nabla\hat F^\mathrm{sm}_{m,V}(z)+\zeta$ are normal random variables conditioned on $V$ and $x$.
    \item $b_y = 2\alpha$, $b_z = 2\beta$.
    \end{itemize}
\end{itemize}
As a result, we have
\begin{equation*}
\begin{split}
\Delta_T^*(\mathcal{F}^0_\mathrm{CSO}, \Phi_\mathrm{m}, \mathcal{A}) 
\geq & \alpha(1-\sqrt{\frac{T\alpha^2}{\sigma^2}}) +\frac{1}{\sqrt{m}}\sqrt{\frac{2}{\pi}}(1-\sqrt{\frac{T\beta^2}{\sigma^2}}).
\end{split}
\end{equation*}
Thus there exists a hard instance with $\alpha = \frac{1}{2}\sqrt{\frac{\sigma^2}{T}}$and $\beta = \frac{1}{2}\sqrt{\frac{\sigma^2}{T}}$, such that 
\begin{equation}
\Delta_T^*(\mathcal{A}, \mathcal{F}^0_\mathrm{CSO}, \Phi_\mathrm{m}) 
\geq  
\sqrt{\frac{\sigma^2}{64T}} + \frac{1}{2\sqrt{m}}\sqrt{\frac{2}{\pi}}
\end{equation}
for convex CSO objective with non-smooth outer function.

Letting the right-hand side of each result be greater or equal to $ \eps $, we have the corresponding sample complexity for each case.  

For the nonconvex CSO problems, we construct a hard instance such that $F^c(y)$ is nonconvex smooth and $F^\mathrm{sm}(z)$ is $1$-strongly convex with Lipschitz continuous or Lipschitz smooth outer function.
Since $F^\mathrm{sm}(z)$ is $1$-strongly convex, it holds
$$
\|\nabla F^\mathrm{sm}(z)\|_2^2\geq 2(F^\mathrm{sm}(z)-\inf_z F^\mathrm{sm}(z)).
$$
We further have
\begin{equation*}
\label{eq:nonconvex_decomposition}
\begin{split}
& \Delta_T^{*g}({\mathcal{A}},\mathcal{F}^-_\mathrm{CSO},\Phi_m,)\\
=&  
\inf_{A \in { }{\mathcal{A}}} \sup_{\phi\in\Phi_m}\sup_{F \in \mathcal{F}^-_\mathrm{CSO}}  \mathbb{E}\|\nabla F\big(x^{A}_{T}(\phi)\big)\|_2^2\\
\geq & 
\inf_{A \in { }{\mathcal{A}}} \sup_{\phi\in\Phi}\sup_{F^c \in \mathcal{F}^-}  \mathbb{E}\|\nabla F^\mathrm{c}\big(y^{A}_{T}(\phi)\big)\|_2^2+ \inf_{A \in { }{\mathcal{A}}} \sup_{\phi\in\Phi_m}\sup_{F^{sm} \in \mathcal{F}^+_\mathrm{CSO}} \mathbb{E}\|\nabla F^\mathrm{sm}\big(z^{A}_{T}(\phi)\big)\|_2^2\\
\geq &
\sup_{\phi\in\Phi}
\sup_{P_F\in\mathcal{P}\{\mathcal{F}^-\}}
\inf_{A\in { }{\mathcal{A}}} \mathbb{E}\|\nabla F\big(y^{A}_{T}(\phi)\big)\|_2^2
+ 
\inf_{A \in { }{\mathcal{A}}} \sup_{\phi\in\Phi_m}\sup_{F^{sm} \in \mathcal{F}^+_\mathrm{CSO}} \mathbb{E}\|\nabla F^\mathrm{sm}\big(z^{A}_{T}(\phi)\big)\|_2^2\\
\geq & 
\sup_{\phi\in\Phi}
\sup_{P_F\in\mathcal{P}\{\mathcal{F}^-\}}
\inf_{A\in { }{\mathcal{A}}} \mathbb{E}\|\nabla F\big(y^{A}_{T}(\phi)\big)\|_2^2
+
\inf_{A \in { }{\mathcal{A}}} \sup_{\phi\in\Phi_m}\sup_{F^{sm} \in \mathcal{F}^+_\mathrm{CSO}} 2\mathbb{E}(F^\mathrm{sm}\big(z^{A}_{T}(\phi)\big)-\inf_z F^\mathrm{sm}(z)),
\end{split}
\end{equation*}
where $\Phi\subset\Phi_\mathrm{m}$ is the oracle class such that $\EE h(x,\zeta)=F(x)$, the first inequality holds by definition, the second inequality uses the fact from \citet{braun2017lower} that
    \begin{equation}
\inf_{A \in \mathcal{A}} \sup_{\phi\in\Phi}\sup_{F \in \mathcal{F}}   \Delta_T^{g}(A,F,\phi) \geq \underset{\phi\in\Phi}{\sup}\ 
\underset{F\in\mathcal{P}\{\mathcal{F}\}}{\sup}\
\underset{A\in{ \mathcal{A}}}{\inf}
\Delta_T^{g}(A,F,\phi),    
\end{equation}
where $\mathcal{P}(\mathcal{F})$ is the set of all distributions over $\mathcal{F}$, the third inequality holds  by strong convexity.  

Note that the second term on the right hand side of the last inequality is exactly twice the minimax error for strongly convex CSO objectives.
Thus we could use the hard instance construction earlier on the strongly convex $\hat F_{m,V}^\mathrm{sm}$ to lower bound it. When the oracle has access to an $\EE \hat F(x;\xi,\{\eta_j\}_{j=1}^m)$ with smooth outer function $f_\xi$, 
$$
\inf_{A \in { }{\mathcal{A}}} \sup_{\phi\in\Phi_m}\sup_{F^{sm} \in \mathcal{F}^+_\mathrm{CSO}} 2\mathbb{E}(F^\mathrm{sm}\big(z^{A}_{T}(\phi)\big)-\inf_z F^\mathrm{sm}(z^*))\geq \frac{1}{2m}.
$$
When the oracle has access to an $\EE \hat F(x;\xi,\{\eta_j\}_{j=1}^m)$ with nonsmooth outer function $f_\xi$, 
$$
\inf_{A \in { }{\mathcal{A}}} \sup_{\phi\in\Phi_m}\sup_{F^{sm} \in \mathcal{F}^+_\mathrm{CSO}} 2\mathbb{E}(F^\mathrm{sm}\big(z^{A}_{T}(\phi)\big)-\inf_z F^\mathrm{sm}(z))\geq \frac{1}{2m}+\sqrt{\frac{2}{m\pi}}.
$$

As for the first term on the right-hand side, we directly use the results from \citet{arjevani2019lower} to lower bound it. \citet{arjevani2019lower} says that for $\mathcal{F}^-_S$, the class of nonconvex smooth functions, it holds for any $\eps>0$ that
$$
\sup_{\phi\in\Phi}
\sup_{P_F\in\mathcal{P}\{\mathcal{F}^-_S\}}
\inf_{A\in { }{\mathcal{A}}} \mathbb{E}\|\nabla F\big(y^{A}_{t}(\phi)\big)\|_2^2 \geq \eps^2,
$$
for any $t\leq t_\mathrm{max}=\cO( \sigma^2\eps^{-4})$ and 
$$
\sup_{\phi\in\Phi^c}
\sup_{P_F\in\mathcal{P}\{\mathcal{F}^-_S\}}
\inf_{A\in { }{\mathcal{A}}} \mathbb{E}\|\nabla F\big(y^{A}_{t}(\phi)\big)\|_2^2 \geq \eps^2,
$$
for $t\leq t_\mathrm{max}^c=\cO( \sigma^2\eps^{-3})$ where $\Phi^c$ denote an stochastic first-order unbiased oracle class such that any oracle in this class has  Lipschitz continuous gradient estimator. It implies that for any $T$, there exists an $\eps>0$ such that  $T=\eps^{-4}\leq t_\mathrm{max}$ and
$$
\sup_{\phi\in\Phi}
\sup_{P_F\in\mathcal{P}\{\mathcal{F}^-_S\}}
\inf_{A\in { }{\mathcal{A}}} \mathbb{E}\|\nabla F\big(y^{A}_{T}(\phi)\big)\|_2^2 \geq \eps^2.
$$
It further implies that
$$
\sup_{\phi\in\Phi}
\sup_{P_F\in\mathcal{P}\{\mathcal{F}^-\}}
\inf_{A\in { }{\mathcal{A}}}
\mathbb{E}\|\nabla F\big(y^{A}_{T}(\phi)\big)\|_2^2
\geq 
\sup_{\phi\in\Phi}
\sup_{P_F\in\mathcal{P}\{\mathcal{F}^-_S\}}
\inf_{A\in { }{\mathcal{A}}}
\mathbb{E}\|\nabla F\big(y^{A}_{T}(\phi)\big)\|_2^2 \geq \cO(\sigma T^{-1/2}),
$$
and
$$
\sup_{\phi\in\Phi^c}
\sup_{P_F\in\mathcal{P}\{\mathcal{F}^-\}}
\inf_{A\in { }{\mathcal{A}}}
\mathbb{E}\|\nabla F\big(y^{A}_{T}(\phi)\big)\|_2^2
\geq 
\sup_{\phi\in\Phi^c}
\sup_{P_F\in\mathcal{P}\{\mathcal{F}^-_S\}}
\inf_{A\in { }{\mathcal{A}}}
\mathbb{E}\|\nabla F\big(y^{A}_{T}(\phi)\big)\|_2^2 \geq \cO(\sigma T^{-2/3}).
$$

\end{proof}

\section{Experiments}
\label{secapp:numerical}
The platform used for the experiments is Intel Core i9-7940X CPU @ 3.10GHz, 32GB RAM, 64-bit Ubuntu 18.04.3 LTS.

\subsection{Invariant Logistic Regression}
\label{secapp:exp_log}
We generate a synthetic dataset  with $d=10$, $ a\sim\mathcal{N}(0; \sigma_1^2I_d)$ with $ \sigma_1^2=1$ , $\eta|\xi \sim  \mathcal{N}(a; \sigma_2^2I_d)$.  Three different variances of $\eta|\xi$: $ \sigma_2^2 \in\{1,10,100\}$ are considered, corresponding to different noise ratios. At each iteration, we use a fixed mini-batch size $m_t=m$, namely $m$ samples of $\eta|\xi$ are generated for a given feature label pair $\xi_i=(a_i,b_i)$. We fine-tune the stepsizes for BSGD using grid search. 

Table \ref{tab:BSGD_SAA_Logistic_FULL} compares the performance achieved by BSGD and SAA under the metric of optimality gap, $F(x)-F^*$. Since we do not have direct access to the function value. We estimate the objective with $50000$ outer samples and calculate the true conditional expectation. The empirical risk minimization constructed by SAA is solved using CVXPY.

\begin{center}
\vskip -0.1in
\renewcommand\arraystretch{1.1}
\begin{table*}[htbp]
    \caption{Comparison of BSGD and SAA in Invariant Logistic Regression}
    \small
        \centering
        \begin{tabular}{cccccccc}
            \hline
            \multicolumn{7}{c}{$\sigma_1=1, \sigma_2=1$} \\
            \cline{1-7}
            \multirow{2}{*}{$m$} & \multicolumn{2}{c}{$Q=10^5$} & \multicolumn{2}{c}{$Q=5\times10^5$} & \multicolumn{2}{c}{$Q=10^6$} \\
            \cline{2-7}
            & \textit{Mean} & \textit{Dev} & \textit{Mean} & \textit{Dev} &\textit{Mean} & \textit{Dev} \\
            \hline
            1 & \textbf{9.28e-04} & 1.95e-04 & 6.23e-04 & 8.18e-05 & 5.81e-04 & 4.00e-05 \\
            5 & 1.04e-03 & 3.06e-04 & \textbf{2.08e-04} & 6.54e-05 & \textbf{1.77e-04} & 4.70e-05 \\
            10 & 1.22e-03 & 2.15e-04 & 3.69e-04 & 8.14e-05 & 2.91e-04 & 4.91e-05 \\
            20 & 1.46e-03 & 8.94e-04 & 3.22e-04 & 1.54e-04 & 1.66e-04 & 6.44e-05 \\
            50 & 1.53e-02 & 3.47e-03 & 8.82e-04 & 3.56e-04 & 3.94e-04 & 1.61e-04 \\
            100 & 3.40e-02 & 8.58e-03 & 1.94e-03 & 6.48e-04 & 9.27e-04 & 3.45e-04 \\
            \hline
            \textbf{SAA (m=100)} & 2.55e-03 & 9.34e-04 & 8.95e-04 & 3.78e-04 & 5.56e-04 & 2.81e-04 \\
            \hline
        \end{tabular}
        \begin{tabular}{cccccccc}
            \hline
            \multicolumn{7}{c}{$\sigma_1=1, \sigma_2=10$} \\
            \cline{1-7}
            \multirow{2}{*}{$m$} & \multicolumn{2}{c}{$Q=10^5$} & \multicolumn{2}{c}{$Q=5\times 10^5$} & \multicolumn{2}{c}{$Q=10^6$} \\
            \cline{2-7}
            & \textit{Mean} & \textit{Dev} & \textit{Mean} & \textit{Dev} &\textit{Mean} & \textit{Dev} \\
            \hline
             1 & 2.47e-03 & 1.12e-03 & 1.02e-03 & 2.83e-04 & 8.16e-04 & 1.38e-04 \\
             5 & \textbf{2.21e-03} & 9.22e-04 & \textbf{5.53e-04} & 1.30e-04 & \textbf{3.26e-04} & 1.15e-04 \\
             10 & 2.32e-03 & 5.29e-04 & 7.22e-04 & 2.55e-04 & 5.32e-04 & 1.72e-04 \\
             20 & 3.57e-03 & 7.88e-04 & 7.37e-04 & 3.25e-04 & 3.99e-04 & 1.37e-04 \\
             50 & 7.87e-03 & 2.96e-03 & 1.42e-03 & 7.57e-04 & 7.25e-04 & 3.65e-04 \\
             100 & 1.91e-02 & 6.46e-03 & 2.23e-03 & 1.01e-03 & 8.90e-04 & 4.83e-04 \\
            \hline
            \textbf{SAA (m=464)} & 8.69e-03 & 2.74e-03 & 3.70e-03 & 1.07e-03 & 2.14e-03 & 8.45e-04 \\
            \hline
        \end{tabular}
        \begin{tabular}{cccccccc}
            \hline
            \multicolumn{7}{c}{$\sigma_1=1, \sigma_2=100$} \\
            \cline{1-7}
            \multirow{2}{*}{$m$} & \multicolumn{2}{c}{$Q=10^5$} & \multicolumn{2}{c}{$Q=5\times 10^5$} & \multicolumn{2}{c}{$Q=10^6$} \\
            \cline{2-7}
            & \textit{Mean} & \textit{Dev} & \textit{Mean} & \textit{Dev} &\textit{Mean} & \textit{Dev} \\
            \hline
            1 & 7.32e-02 & 7.94e-03 & 6.82e-02 & 2.41e-03 & 6.69e-02 & 1.09e-03 \\
            5 & 1.53e-02 & 4.54e-03 & 3.30e-03 & 1.12e-03 & 1.61e-03 & 7.50e-04 \\
            10 & \textbf{1.46e-02} & 3.80e-03 & 3.28e-03 & 1.24e-03 & 1.70e-03 & 5.82e-04 \\
            20 & 1.73e-02 & 8.95e-03 & \textbf{3.19e-03} & 1.18e-03 & 1.52e-03 & 5.60e-04 \\
            50 & 1.47e-02 & 5.15e-03 & 3.36e-03 & 1.27e-03 & \textbf{1.50e-03} & 6.97e-04 \\
            100 & 3.20e-02 & 8.07e-03 & 5.81e-03 & 2.44e-03 & 3.39e-03 & 1.30e-03 \\
            \hline
            \textbf{SAA (m=1000)} & 4.33e-02 & 1.19e-03 & 1.50e-02 & 8.00e-04 & 1.12e-02 & 6.42e-04 \\
            \hline
        \end{tabular}
    \label{tab:BSGD_SAA_Logistic_FULL}
    \end{table*}
\end{center}

\subsection{MAML}
\label{secapp:exp_MAML}
We use the objective value as the measurement. Since the objective is analytically intractable, we evaluate the MAML objective  via empirical objective obtained by empirical risk minimization:
\begin{equation}
\label{eq:empirical_maml}
    \hat{F}(w)=\frac{1}{\hat{T}}\sum_{i=1}^{\hat{T}}\frac{1}{\hat{N}}\sum_{n=1}^{\hat{N}} l_i\Big(w-\alpha\cdot\frac{1}{\hat{M}}\sum_{m=1}^{\hat{M}}\nabla_w l_i(w,D_{support}^{i,m}); D_{query}^{i,n}\Big),
\end{equation}
where the three sample sizes $ \hat{T} $, $ \hat{N} $ and $ \hat{M} $ are set to be 100. when computing the approximate loss function value, the sample tasks/data are selected randomly.

\begin{figure}
    \centering
\includegraphics[width=.35\textwidth,trim=10 10 80 40,clip]{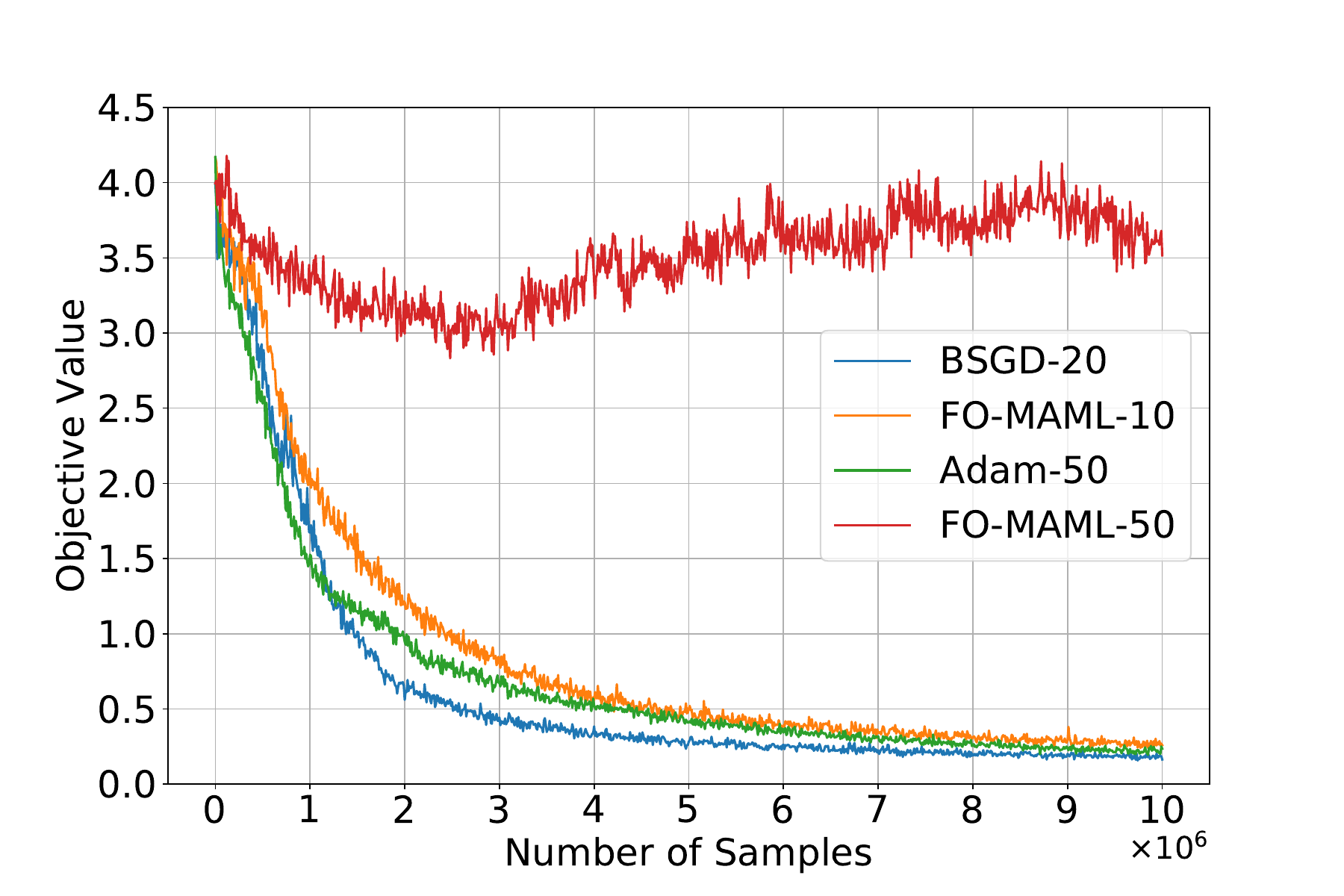}
    \caption{FO-MAML may not converge}
    \label{fig:FO-MAML}
\end{figure}

Figure \ref{fig:FO-MAML}  shows that the widely used first-order MAML~\citep{finn2017model}, which ignores the Hessian information when constructing the gradient estimator, may not converge.  The number after each method denotes the inner mini-batch size. It compares the convergences of the widely used First-order MAML(FO-MAML)\citep{finn2017model}, BSGD, and  Adam, each under the best-tuned inner batch size. BSGD achieves the least error among the three methods with a proper inner batch size of $20$.  Adam requires a larger inner batch size to achieve its best performance, which is less practical as some tasks only have a few or even one sample. 

Table \ref{tab:BSGD_MAML_Value_FULL} summarizes the detailed experimental results of BSGD, FO-MAML, Adam, and BSpiderBoost with different inner mini-batch sizes. The total sample size is $Q=10^7$. The stepsizes for BSGD, FO-MAML, and BSpiderBoost are fine-tuned. Specifically for BSpiderBoost we set $ (N_1=10, N_2=1, q=10) $.  For each inner mini-batch size, we run each algorithm for $10$ times and then calculate the mean and the standard deviation of the output objectives of all trials. The best performance result for each algorithm is highlighted using bold font.

\begin{center}
\renewcommand\arraystretch{1.1}
    \begin{table*}[htbp]
        \caption{Comparison of convergence results of BSGD, FO-MAML and Adam in MAML problem with different inner mini-batch sizes. }
        \small
        \centering
        \begin{tabular}{cccccccc}
            \hline
            \multicolumn{8}{c}{$\alpha=0.01$} \\
            \cline{1-8}
            \multirow{2}{*}{\textbf{Method}} 
            & \multirow{2}{*}{$m$} & \multicolumn{2}{c}{$Q=10^5$} & \multicolumn{2}{c}{$Q=10^6$} & \multicolumn{2}{c}{$Q=10^7$} \\
            \cline{3-8}
             & & \textit{Mean} & \textit{Dev} & \textit{Mean} & \textit{Dev} &\textit{Mean} & \textit{Dev} \\
            \hline
            \multirow{5}{*}{{\textbf{BSGD}}} 
             & 5 & 3.46e+00 & 1.81e-01 & 1.28e+00 & 1.72e-01 & 5.07e-01 & 1.02e-01 \\
             & 10 & \textbf{3.40e+00} & 2.43e-01 & \textbf{1.20e+00} & 2.68e-01 & 2.68e-01 & 8.09e-02 \\
             & 20 & 3.57e+00 & 3.18e-01 & 1.67e+00 & 6.87e-01 & \textbf{1.63e-01} & 5.82e-02 \\
             & 50 & 3.44e+00 & 2.07e-01 & 2.51e+00 & 6.12e-01 & 2.51e-01 & 5.21e-02 \\
             & 100 & 3.81e+00 & 3.66e-01 & 3.23e+00 & 2.99e-01 & 3.60e-01 & 9.05e-02 \\
            \hline
            \multirow{5}{*}{{\textbf{FO-MAML}}}
             & 5 & 3.89e+00 & 3.46e-01 & 3.21e+00 & 2.49e-01 & 8.48e-01 & 1.59e-01 \\
             & 10 & \textbf{3.72e+00} & 3.18e-01 & \textbf{2.07e+00} & 4.94e-01 & \textbf{2.58e-01} & 4.02e-02 \\
             & 20 & 4.03e+00 & 3.32e-01 & 3.15e+00 & 1.69e-01 & 1.82e+00 & 7.12e-01 \\
             & 50 & 3.90e+00 & 3.84e-01 & 3.26e+00 & 3.24e-01 & 3.52e+00 & 5.49e-01 \\
             & 100 & 3.80e+00 & 3.62e-01 & 3.48e+00 & 2.51e-01 & 4.09e+00 & 5.10e-01 \\
            \hline
            \multirow{5}{*}{{\textbf{Adam}}}
             & 5 & \textbf{2.95e+00} & 5.90e-01 & 1.45e+00 & 4.15e-01 & 1.04e+00 & 3.83e-01\\
             & 10 & 3.03e+00 & 4.26e-01 & 1.34e+00 & 5.61e-01 & 6.09e-01 & 7.59e-01 \\
             & 20 & 3.47e+00 & 3.01e-01 & \textbf{1.11e+00} & 3.63e-01 & 2.82e-01 & 8.85e-02 \\
             & 50 & 3.36e+00 & 3.43e-01 & 1.53e+00 & 5.20e-01 & \textbf{2.35e-01} & 8.32e-02 \\
             & 100 & 3.60e+00 & 2.86e-01 & 2.52e+00 & 5.28e-01 & 3.92e-01 & 2.20e-01 \\
            \hline
            \multirow{4}{*}{\textbf{BSpiderBoost}}
			& 10 & 2.47e+00 & 1.02e+00 & 2.65e+00 & 1.08e+00 & 2.56e+00 & 1.04e+00 \\
			& 20 & 6.81e-01 & 6.50e-01 & \textbf{1.95e-01} & 4.49e-02 & \textbf{1.43e-01} & 2.57e-02 \\
			& 50 & \textbf{6.15e-01} & 1.93e-01 & 2.54e-01 & 6.38e-02 & 2.10e-01 & 3.56e-02 \\
			& 100 & 3.21e+00 & 1.12e+00 & 2.76e+00 & 1.45e+00 & 2.98e+00 & 1.55e+00 \\
			\hline
        \end{tabular}
    \label{tab:BSGD_MAML_Value_FULL}
    \end{table*}
\end{center}

\end{appendix}
\end{document}